\documentclass[12pt]{article}
\usepackage{color}
\usepackage{euscript,amssymb,amsfonts,amsmath,amsthm}
\usepackage{mathrsfs}
\usepackage{mathtext}
\usepackage{mathrsfs}
\usepackage[cp1251]{inputenc}
\usepackage[T2A]{fontenc}
\usepackage[english,russian]{babel}

\renewcommand{\P}{{\sf P}}
\newcommand{\E}{{\sf E}\,}
\newcommand{\p}{{\sf P}}
\newcommand{\e}{{\sf E}\,}
\newcommand{\D}{{\sf D}\,}
\newcommand{\pto}{\stackrel{P}{\longrightarrow}}
\newcommand{\eqd}{\stackrel{d}{=}}
\renewcommand{\r}{\mathbb{R}}

\newcommand{\il}[2]{\int\limits_{#1}^{#2}}
\newcommand{\eps}{\varepsilon}

\renewcommand{\kappa}{\varkappa}
\renewcommand{\le}{\leqslant}
\renewcommand{\ge}{\geqslant}

\newcommand{\R}{\mathbb R}

\newcommand{\I}{{\bf1}} %индикатор

\newcommand{\F}{\mathcal F}
\newcommand{\ud}{\rho(F_n,\Phi)} %uniform distance

\newcommand{\la}{\lambda}
\newcommand{\bet}{\beta^3}

\newcommand{\sign}{\mathrm{sign}}

\newcommand{\si}{\sigma}

\theoremstyle{remark}

\renewcommand{\proof}{{P r o o f}}

\bibstyle{plain}

\title{An improvement of the Berry--Esseen inequality \\ with
applications to Poisson and mixed Poisson random
sums\thanks{Research supported by the Russian Foundation for Basic
Research, projects 08-01-00563, 08-01-00567, 08-07-00152 and
09-07-12032-ofi-m, and also by the Agency for Education of Russian Federation, state contracts P-1181 and P-958.}}
\author{Victor Korolev\footnote{Department of Mathematical Statistics, Faculty of Computational
Mathematics and Cybernetics, Moscow State University; Institute
for Informatics Problems, Russian Academy of Sciences,
vkorolev@comtv.ru} and Irina Shevtsova\footnote{Department of
Mathematical Statistics, Faculty of Computational Mathematics and
Cybernetics, Moscow State University, ishevtsova@cs.msu.su} }

\date{15 December 2009}

\begin{document}

\maketitle

{\small

\centerline{\bf Abstract.}

\smallskip

By a modification of the method that was applied in (Korolev and
Shevtsova, 2009), here the inequalities
$$
\ud\le\frac{0.335789(\beta^3+0.425)}{\sqrt{n}}
$$
and
$$
\ud\le \frac{0.3051(\beta^3+1)}{\sqrt{n}}
$$
are proved for the uniform distance $\rho(F_n,\Phi)$ between the
standard normal distribution function $\Phi$ and the distribution
function $F_n$ of the normalized sum of an arbitrary number
$n\ge1$ of independent identically distributed random variables
with zero mean, unit variance and finite third absolute moment
$\beta^3$. The first of these inequalities sharpens the best known
version of the classical Berry--Esseen inequality since
$0.335789(\beta^3+0.425)\le0.335789(1+0.425)\beta^3<0.4785\beta^3$
by virtue of the condition $\beta^3\ge1$, and $0.4785$ is the best
known upper estimate of the absolute constant in the classical
Berry--Esseen inequality. The second inequality is applied to
lowering the upper estimate of the absolute constant in the analog
of the Berry--Esseen inequality for Poisson random sums to
$0.3051$ which is strictly less than the least possible
value of the absolute constant in the classical Berry--Esseen
inequality. As a corollary, the estimates of the rate of
convergence in limit theorems for compound mixed Poisson
distributions are refined.

\smallskip


}

\section{Introduction}

By $\F_3$ we will denote the set of distribution functions with
zero first moment, unit second moment and finite third absolute
moment $\beta^3$. Let $X_1,X_2,\ldots$~be independent random
variables with common distribution function $F\in\F_3$ defined on
a probability space $(\Omega,\mathcal{A},\P)$. Denote
$$
F_n(x)=F^{*n}(x\sqrt{n})=
\P\left(\frac{X_1+\ldots+X_n}{\sqrt{n}}<x\right),
$$
$$
\Phi(x)=\int_{-\infty}^x\phi(t)dt,\quad
\phi(x)=\frac1{\sqrt{2\pi}}e^{-x^2/2},\quad x\in\R.
$$
The classical Berry--Esseen theorem states that there exists a
finite positive absolute constant $C_0$ which guarantees the
validity of the inequality
$$
\ud\equiv\sup_x|F_n(x)-\Phi(x)|\le
C_0\frac{\beta^3}{\sqrt{n}}\eqno(1)
$$
for all $n\ge1$ and any $F\in\F_3$ (Berry, 1941), (Esseen, 1942).
The problem of establishing the best value of the constant $C_0$
in inequality (1) is very important from the point of view of
practical estimation of the accuracy of the normal approximation
for the distribution functions of random variables which may be
assumed to have the structure of a sum of independent random
summands.

This problem has a long history and is very rich in deep and
interesting results. Upper estimates for $C_0$ were considered in
very many papers. Here we will not repeat a detailed history of
the efforts to lower the upper estimates of $C_0$ from the
original works of A. Berry (Berry, 1941) and C.-G. Esseen (Esseen,
1942) to the papers of I. S. Shiganov (Shiganov, 1982), (Shiganov,
1986) presented in (Korolev and Shevtsova, 2009). We will
restrict ourselves only to an outline of the recent history of the
subject.

After some lull that lasted more than twenty years, recently the
interest to the problem of improving the Berry--Esseen inequality
rose again and resulted in very interesting and in some sense
path-clearing works. In 2006 I. G. Shevtsova improved Shiganov's
upper estimate by approximately $0.06$ and obtained the estimate
$C_0\le0.7056$ (Shevtsova, 2006). In 2008 she sharpened this
estimate to $C_0\le0.7005$ (Shevtsova, 2008). In 2009 the
competition for improving the constant became especially keen. On
8 June, 2009 I. S. Tyurin submitted his paper (Tyurin, 2009a) to
the <<Theory of Probability and Its Applications>>. That paper,
along with other results, contained the estimate $C_0\le0.5894$.
Two days later the summary of those results was submitted to
<<Doklady Akademii Nauk>> (translated into English as <<Doklady
Mathematics>>) (Tyurin, 2009b). Independently, on 14 September,
2009 V. Yu. Korolev and I. G. Shevtsova submitted their paper
(Korolev and Shevtsova, 2009) to the <<Theory of Probability and
Its Applications>>. In that paper the inequality
$$
\ud\le \frac{0.34445(\beta^3+0.489)}{\sqrt{n}}, \ \ \
n\ge1,\eqno(2)
$$
was proved which holds for any distribution $F\in\F_3$ yielding
the estimate $C_0\le0.5129$ by virtue of the condition
$\beta^3\ge1$. Finally, on 17 November, 2009 the paper (Tyurin,
2009c) was submitted to the <<Russian Mathematical Surveys>> {(its
English version (Tyurin, 2009d) appeared on 3 December, 2009 on
arXiv:0912.0726v1)}. In this paper the estimate $C_0\le0.4785$ is
proved. So, the best known upper estimate of the absolute constant
$C_0$ in the classical Berry--Esseen inequality (1) is
$C_0\le0.4785$ (Tyurin, 2009c).

On the other hand, in 1956 C.-G. Esseen showed that $C_0\ge C_E$
where
$$
C_{E}=\frac{\sqrt{10}+3}{6\sqrt{2\pi}}=0.409732...
$$
(Esseen, 1956). In 1967 V. M. Zolotarev put forward the hypothesis
that in (1) $C_0=C_E$ (Zolotarev, 1967a), (Zolotarev, 1967b).
However, up till now this hypothesis has been neither proved nor
rejected.

To prove (2) we used an observation that from inequality (1) it
obviously follows that for any $k\ge0$ there exists a finite
positive absolute constant $C_k$ which guarantees the validity of
the inequality
$$
\ud\le C_k\frac{\beta^3+k}{\sqrt{n}}\eqno(3)
$$
for all $n\ge1$ and $F\in\F_3$ (for example, inequality (3)
trivially holds with $C_k=C_0$).

Following the lines of the reasoning we used in (Korolev and
Shevtsova, 2009) to prove (2), with the only change in the way of
estimation of the difference between characteristic functions in
the neighborhood of zero (see lemma 2 below), in this paper we
will demonstrate a special method of numerical estimation of $C_k$
in (3). This method yields two special values of $k$: $k=k_0$ and
$k=1$. The first value, $k_0$, minimizes the upper estimate of
$C_k(1+k)$ yielding the best (within the method under
consideration) upper estimate of $C_0$ in (1) since
$$
C_0\le\min_{k\ge0}C_k(1+k)
$$
by virue of the condition $\beta^3\ge1$. At the same time the
second value, $k=1$, minimizes $C_k$ in (3). As we will see, $k=1$
plays the main role in improving the absolute constant in the
analog of the Berry--Esseen inequality for Poisson and mixed
Poisson random sums.

Inequality (3) with $k=k_0$ and $k=1$ is an improvement of the
inequality
$$
\ud\le 0.3450\frac{\beta^3+1}{\sqrt{n}}
$$
we proved in (Korolev and Shevtsova, 2010a). In (Korolev and
Shevtsova, 2010b) this inequality was applied to sharpening the
analog of the Berry--Esseen inequality for Poisson random sums and
it was for the first time demonstrated that the absolute constant
in this analog can be made strictly less than that in the
classical Berry--Esseen inequality.

In the papers (Shevtsova, 2010a) and (Korolev and Shevtsova,
2010a) it was shown that the constant $C_k$ in (3) cannot be made
less than the so-called lower asymptotically exact constant in the
central limit theorem, that is,
$$
C_k\ge \frac{2}{3\sqrt{2\pi}}=0.2659...,
$$
so that the gaps between the least possible value of the constant
$C_k$ and its upper estimates given in theorems 1 and 2 below are
rather small and do not exceed 0.07 and 0.035, respectively, which
is important from the point of view of practical applications of
inequalities (6) and (7).

Our investigations were to a great extent motivated by a series of
results of H{$^{^{\circ}}\hspace{-2mm}$a}kan Prawitz and Vladimir
Zolotarev outlined below.

First, since estimates of the accuracy of the normal approximation
for distributions of sums of independent random variables are
traditionally constructed with the use of the so-called smoothing
inequalities which estimate the (uniform) distance between the
pre-limit distribution function of the standardized sum of
independent random variables and the limit standard normal
distribution function through some integral of the (weighted)
absolute value of the difference between the corresponding
characteristic functions, the shape of the dependence of the final
estimate on the moments of summands is fully determined by the
shape of dependence of the majorant of characteristic functions on
these moments. In (Prawitz, 1973) the following result was
presented. Let $f(t)$ be the characteristic function corresponding
to the distribution function $F\in\F_3$. Denote
$$
\kappa=\sup_{x>0}\frac{\left|\cos
x-1+x^2/2\right|}{x^3}=0.09916191...
$$
and let $\theta_0=3.99589567...$~be the unique root of
the equation
$$
3(1-\cos\theta)-\theta\sin\theta-\theta^2/2=0,
$$
lying in the interval $(\pi,2\pi)$. Then
$$
|f(t)|\le \left\{%
\begin{array}{ll}
\displaystyle 1-\frac{t^2}2+\kappa\left(\beta_3+
 1\right)|t|^3, &\displaystyle
 |t|\le\frac{\theta_0}{(\beta_3+1)}, \\
\displaystyle 1-\frac{1-\cos\big((\beta_3+1) t\big)}
 {(\beta_3+1)^2},
 & \theta_0\le (\beta_3+1)|t|  \le2\pi, \\
    1, &\displaystyle |t|\ge\frac{2\pi}{(\beta_3+1)}.
\end{array}%
\right.
$$
As is easily seen, the majorant for $|f(t)|$ established by this
inequality depends on $\beta_3$ through the function
$\psi(\beta_3)=\beta_3+1$. This is the first hint at that the
final estimate for $\rho(F_n,\Phi)$ should also depend on
$\beta_3$ through the function $\psi(\beta_3)=\beta_3+1$.

Second, in (Prawitz, 1975b) H. Prawitz announced an inequality
with unusual structure
$$
\rho(F_n,\Phi)\le
\frac{2}{3\sqrt{2\pi}}\cdot\frac{\beta_3}{\sqrt{n-1}}+
\frac{1}{2\sqrt{2\pi(n-1)}}+
\frac{c_1(\beta_3)^2+c_2\beta_3+c_3}{n-1},\quad n\ge2,\ F\in\F_3,
\eqno(4)
$$
where $c_1$, $c_2$ and $c_3$~are some finite positive constants.
In the same paper he suggested that the coefficient
$$
\frac{2}{3\sqrt{2\pi}}=0.2659...
$$
at~$\beta_3/\sqrt{n-1}$ cannot be made smaller. Probably, H.
Prawitz intended to publish the strict proof of ~(4) in the second
part of his work which, unfortunately, for some reasons remained
unpublished (the title of (Prawitz, 1975b) contains the Roman
number I indicating the assumed continuation).

This Prawitz' inequality (4) seemed to have bepuzzled some
specialists in limit theorems of probability theory. In
particular, it was bypassed in the well-known books (Petrov,
1987), (Zolotarev, 1997) (in both of these books there is even no
reference to any of Prawitz' works). Only in the book (Petrov,
1995) there appears a reference to the paper (Prawitz, 1975a)
dealing with some estimates for characteristic functions, but the
paper (Prawitz, 1975b) containing inequality (4) is again ignored.
In Mathematical Reviews (Dunnage, 1977) there is only a fuzzy
remark concerning <<some improvements for identically distributed
summands>>. Probably, this attitude of some specialists to
inequality~(4) is caused by that at first sight this inequality
contradicts the Esseen's result that $C_0\ge C_E$ cited above,
since
$$
\frac2{3\sqrt{2\pi}}< \frac{\sqrt{10}+3}{6\sqrt{2\pi}}.
$$
However, a thorough analysis of the published part of Prawitz'
work convinces that inequality (4) is valid. A strict proof of a
similar inequality for not necessarily identically distributed
summands with the third term being
$O\big((\beta_3/\sqrt{n})^{5/3}\big)$ was given by V.~Bentkus
(Bentkus, 1991), (Bentkus, 1994) (for identically distributed
summands, the result of Benkus is slightly worse than (4) where
the third term is $O\big((\beta_3/\sqrt{n})^2\big)$).

Inequality~(4) has a very interesting structure: from the main
term of order $O(n^{-1/2})$ of the estimate of the accuracy of the
normal approximation a summand of the form $1/\sqrt{n}$ is
separated. This summand may be considerably less than the Lyapunov
fraction $\beta_3/\sqrt{n}$. Moreover, in the double array scheme
it may happen so that even if the Lyapunov condition
$\beta_3/\sqrt{n}\to0$ holds, the quantity $\beta_3=\beta_3(n)$
may infinitely increase as $n\to\infty$ so that the summand of the
form $n^{-1/2}$ is infinitesimal with a higher order of smallness
than the Lyapunov fraction $\beta_3(n)/\sqrt{n}$. Thus, inequality
(4) is the second hint at that in a reasonable estimate of
$\rho(F_n,\Phi)$ depending on $\beta_3$ the term of order
$O(n^{-1/2})$ should be split into two summands of the form
$\beta_3/\sqrt{n}$ and $1/\sqrt{n}$ respectively.

By the way, speaking of the history of inequality (4), it has to
be noted that actually it is a further development of the
inequality
$$
\rho(F_n,\Phi)\le\frac{0.32\beta_3+0.25}{\sqrt{n-2}},\ \ \
n\ge3,\eqno(5)
$$
which holds under the condition $\sqrt{n-1}\ge3.9(\beta_3+1)$. The
proof of (5) was given by H.\,Prawitz in his lecture on 16 June,
1972 at the Summer School of the Swedish Statistical Society in
L{\"o}ttorp (Prawitz, 1972a).

So, the final shape of inequality (3) was prompted by the works of
H. Prawitz mentioned above. As this is so, the main role goes to
the problem of a proper estimation of the constant $C_k$. To solve
this problem we use a method which is a further development of the
ideas of V. Zolotarev presented in (Zolotarev, 1965), (Zolotarev,
1966), (Zolotarev, 1967a) and (Zolotarev, 1967b). This method will
be described in detail below.

The paper is organized as follows. In Section 2 the basic results
are proved. Namely, here we prove inequality (3) with
$k=k_0=0.425$ (theorem 1) and with $k=1$ (theorem 2). In Section 3
theorem 2 is applied to sharpening the analog of the Berry--Esseen
inequality for Poisson random sums. We show that despite a
prevalent opinion that the absolute constant in this inequality
should not be less than the absolute constant in the classical
Berry--Esseen inequality, as a matter of fact this is not so and
the constant in the Berry--Esseen inequality for Poisson random
sums does not exceed 0.3051, which is, as it has been
already mentioned, strictly less than the least possible value
$C_E\approx0.4097$ of the constant $C_0$ in (1). Finally, in
Sections 4 and 5 the result of Section 3 is used for improving the
estimates of the rate of convergence of compound mixed Poisson
distributions with zero and non-zero means to scale and location
mixtures of normal laws, respectively.

\section{The basic results}

\subsection{Formulations and discussion}

Practical calculations show that {\it under the algorithm we use
for the estimation of}  $\,C_k$ (see Section 3) the resulting
majorant of the constant $C_k$ decreases as $k$ increases from $0$
to~$1$. At the same time for $0\le k\le k_0\approx0.425$ the
obtained estimates of $C_k(1+k)$ remain constant, and for $k>k_0$
they begin to increase although in the interval $k_0<k<1$ the
obtained estimate of $C_k$ decreases. Thus, we can present two
computationally optimal values of $k$ in (3): $k_0=0.425$ and
$k_1=1$. The first of them delivers the minimum value to the upper
estimate of $C_k(1+k)$, thus solving the problem of estimation of
$C_0$ in (1), whereas the second, maximin, minimizes the estimate
of $C_k$ in (3).

The use of $k=k_0$ in (3) gives the following result.

\smallskip

{\sc Theorem 1.} {\it For all $n\ge1$ and all distributions with
zero mean, unit variance and finite third absolute moment
$\beta^3$ we have the inequality}
$$
\ud\le \frac{0.335789(\beta^3+0.425)}{\sqrt{n}}.\eqno(6)
$$

\smallskip

{\sc Remark 1}. Under the conditions imposed on the moments of the
random variable $X_1$ we always have $\beta^3\ge1$. Therefore,
$$
0.335789(\beta^3+0.425)\le 0.335789(1+0.425)\beta^3 <
0.4785\beta^3.
$$
Hence, inequality (6) is always sharper than the classical
Berry--Esseen inequality (1) with the best known constant
$C_0=0.4785$ for all possible values of $\beta^3$, although the
same prior information concerning the distribution $F$ is required
for its validity (namely, only the value of the third
absolute moment $\beta^3$).

\smallskip

{\sc Remark 2.} Inequality (6) is an <<unconditional>> variant of
the <<conditional>> Prawitz inequality (5) and is a practically
computable analog of inequality (4) with a slightly (approximately
by 0.07) worse first coefficient and a slightly better
(approximately by 0.05) second coefficient, but without the third
summand that contains unknown constants.

\smallskip

{\sc Remark 3.} Even if the hypothesis of V. M. Zolotarev that
$C_0=C_{E}=0.4097...$ in (1) (see (Zolotarev, 1967a), (Zolotarev,
1967b)) turns out to be true, then, due to that $\beta^3\ge1$,
inequality (6) will be sharper than the classical Berry--Esseen
inequality (1) for $\beta^3\ge 1.93$.

The use of $k=1$ in (3) yields the following result.

\smallskip

{\sc Theorem 2.} {\it For all $n\ge1$ and all distributions with
zero mean, unit variance and finite third absolute moment
$\beta^3$ we have the inequality}
$$
\ud\le \frac{0.3051(\beta^3+1)}{\sqrt{n}}.\eqno(7)
$$

\smallskip

{\sc Remark 4.} Inequality (7) is another <<unconditional>>
variant of the <<conditional>> Prawitz inequality (5). Moreover,
the first coefficient in (7) is less than that in (5) by
approximately $0.02$ whereas the second coefficient in (7) is
greater than that in (5) by approximately $0.05$.

\subsection{Proofs of basic results}

\subsubsection{Auxiliary statements}

As we have already mentioned above, to prove theorem~1 we will
follow the lines of the approach proposed and developed by
V.~M.~Zolotarev in his works (Zolotarev, 1965), (Zolotarev, 1966)
and (Zolotarev, 1967). This approach is based on the application
of smoothing inequalities which make it possible to estimate the
distance between distribution functions via the distances between
the corresponding characteristic functions. Within this approach
the key points are: (i) the choice of a proper smoothing
inequality; (ii) the choice of a proper smoothing kernel in a
smoothing inequality; (iii) the choice of proper estimates for the
distance between characteristic functions; (iv) the choice of a
proper computational optimization procedure.

We will describe these points one after another as they are used
in the proof of theorems~1 and~2. The corresponding
statements will have the form of lemmas.

We begin with the smoothing inequality. In most papers dealing
with the estimation of the constant in the Berry--Esseen
inequality (1) smoothing inequalities of the same type were used.
This type of smoothing inequalities was introduced by V. M.
Zolotarev. In the original paper (Zolotarev, 1965), just as in
similar inequalities in the earlier papers of Berry (Berry, 1941)
and Esseen (Esseen, 1942), the kernel was used which had a
probabilistic sense, that is, which was the probability density of
some symmetric probability distribution. In the paper of Van Beek
(Van Beek, 1972) it was noticed that this condition is not
crucial. Van Beek proposed to use symmetric kernels with
alternating signs. Concurrently with (Van Beek, 1972), the paper
of V. Paulauskas (Paulauskas, 1971) was published in which the
original smoothing inequality of Zolotarev was generalized (and
hence, sharpened) to the case of positive non-symmetric kernels.
It is interesting to notice that although in the final part of the
paper of Paulauskas it was noted that the smoothing inequality
proved in that paper was destined, in the first place, for
improving the constant in the Berry--Esseen inequality, as far as
we know, unfortunately no one ever used the Paulauskas inequality
for this purpose. In (Shevtsova, 2009b) a new smoothing inequality
was proved which generalizes (and hence, sharpens) both
Paulauskas' and Van Beek's inequalities to the case of
non-symmetric kernels with alternating signs. However, all these
inequalities yield worse estimates than the Prawitz smoothing
inequality proved in (Prawitz, 1972b).

The characteristic function of the standardized sum
$(X_1+\ldots+X_n)/\sqrt{n}$ will be denoted $f_n(t)$,
$$
f_n(t)=\int\limits_{-\infty}^{\infty}e^{itx}dF_n(x),\ \ \ t\in\r.
$$
Also denote
$$
r_n(t)=|f_n(t)-e^{-t^2/2}|.
$$

\smallskip

{\sc Lemma 1} (Prawitz, 1972b). {\it For an arbitrary distribution
function $F$ and $n\ge1$ for any $0<t_0\le1$ and $T>0$ we have the
inequality
$$ \ud\le 2\int_0^{t_0}|K(t)|r_n(Tt)\,dt+
2\int_{t_0}^{1}|K(t)|\cdot|f_n(Tt)|dt+
$$
$$
+2\int_0^{t_0}\left|K(t)-\frac i{2\pi t}\right|e^{-T^2t^2/2}dt +
\frac1{\pi}\int_{t_0}^\infty e^{-T^2t^2/2}\frac{dt}t,
$$
where
$$ K(t)=\frac12(1-|t|)+\frac i2\left[(1-|t|)\cot\pi
t+\frac{\sign t}\pi\right],\quad -1\le t\le1.\eqno(8)
$$
}

\smallskip

{\sc Remark 5.} In (Vaaler, 1985) a proof of a result similar to
the Prawitz inequality stated by lemma 1 was given by a techniques
different from that used in (Prawitz, 1972b) and it was also
proved that the kernel $K(t)$ defined by (8) is in some sense
optimal.

\smallskip

Now consider the estimates of the characteristic functions
appearing in lemma~1. For $\eps>0$ set
$$
\chi(t,\eps)=\begin{cases}t^2/2-\varkappa\eps|t|^3,&\text{$|t|\le\theta_0/\eps$},\vspace{1mm}\cr
\displaystyle{\frac{1-\cos\eps
t}{\eps^2}},&\text{$\theta_0<\eps|t|\le 2\pi$},\vspace{1mm}\cr
0,&\text{$|t|>2\pi/\eps$},\end{cases}\eqno(9)
$$
where $\theta_0=3.99589567...$ is the unique root of the equation
$$
\theta^2 + 2\theta\sin \theta + 6(\cos \theta - 1)=0,\quad
\pi\le\theta\le2\pi,\eqno(10)
$$
$$
\varkappa \equiv \sup_{x>0}\frac{\left|\,\cos
x-1+x^2/2\right|}{x^3}=\frac{\cos
x-1+x^2/2}{x^3}\left|{\atop{_{x=\theta_0}}}\right.
=0.09916191...\eqno(11)
$$
It can easily be made sure that the function $\chi(t,\eps)$
monotonically decreases in $\eps>0$ for any fixed $t\in\R$.

The Lyapunov fraction will be denoted $\ell=\beta^3/\sqrt{n}$.
In addition, denote
$$
\ell_n=\ell+1/\sqrt{n}.
$$

\smallskip

{\sc Lemma 2}. {\it For any $F\in\F_3$, $n\ge1$ and $t\in\R$ the
following estimates take place:}
$$
|f_n(t)|\le  \Big[1-\frac{2}{n}\chi(t,\ell_n)\Big]^{n/2}
\equiv f_1(t,\ell_n,n),
$$
$$
|f_n(t)|\le \exp\{-\chi(t,\ell_n)\}\equiv  f_2(t,\ell_n),
$$
$$
|f_n(t)|\le \exp\Big\{-\frac{t^2}2+\varkappa\ell_n|t|^3 \Big\}\equiv f_3(t,\ell_n).
$$

\smallskip

{\sc Remark 6.} Apparently, the function $f_1(t,\ell_n,n)$ was
used in the problem of numerical evaluation of the absolute
constants in the estimates of the accuracy of the normal
approximation for the first time in (Korolev and Shevtsova, 2009).
The second and the third estimates presented in lemma 2 are due to
H. Prawitz (Prawitz, 1973), (Prawitz, 1975b).

\smallskip

{\sc Remark 7.} Evidently, $ f_1(t,\eps,n)\le  f_2(t,\eps)$ for
all $n\ge1$, $\eps>0$ and $t\in\R$. Moreover, from the result of
Prawitz (Prawitz, 1973) it follows that $ f_2(t,\eps)\le
f_3(t,\eps)$ for all $\eps>0$ and $t\in\R$, thus the sharpest
estimate for $|f_n(t)|$ is given by $ f_1(t,\ell_n,n)$, while the
estimates $ f_j(t,\ell_n)$, $j=2,3$, possess a useful property of
monotonicity in $\ell_n$ which is very important for the
computational procedure.

\smallskip

{\sc Lemma 3} (Tyurin, 2009a), (Tyurin, 2009c), (Tyurin, 2009d).
{\it For any $F\in\F_3$, $n\ge1$ and $t\in\R$ we have}
$$
r_n(t)\le \ell e^{-t^2/2}\int_0^{|t|} \frac{u^2}2\, e^{u^2/2} \Big|f\Big(\frac{u}{\sqrt{n}}\Big)\Big|^{n-1} du.
$$
\smallskip

The combination of lemmas 2 and 3 allows to obtain an estimate for
the difference of the characteristic functions in the neighborhood
of zero, which is sharper than all the analogous estimates used in
the preceding works:
$$
r_n(t)\le \ell e^{-t^2/2}\int_0^{|t|} \frac{u^2}2\, e^{u^2/2} \Big[1-\frac{2}{n}\chi\Big(u,\ell+
\frac1{\sqrt n}\Big)\Big]^{(n-1)/2} du\equiv r_1(t,\ell,n),\quad t\in\R.
$$
From what was said above it follows that the substitution of the
functions $f_j(t,\ell_n)$, $j=2,3,$ instead of $f_1(t,\ell_n,n)$
into the right-hand side of the last inequality does not make the
resulting estimate less, thus, we obtain two more estimates for
$r_n(t)$ which monotonically increase in $\ell$:
$$
r_n(t)\le \ell e^{-t^2/2}\int_0^{|t|} \frac{u^2}2\, e^{u^2/2}
\exp\Big\{-\frac{n-1}n\cdot\chi\Big(u,\ell+\frac1{\sqrt n}\Big)\Big\} du\equiv r_2(t,\ell,n),
$$
$$
r_n(t)\le \ell e^{-t^2/2}\int_0^{|t|} \frac{u^2}2 \exp\Big\{
\varkappa\ell_n u^3 +
\frac{u^2}{2n}\Big(1-2\varkappa\ell_n u\Big)
\Big\}du\equiv r_3(t,\ell,n),\quad t\in\R,
$$
%$$
%r_n(t)\le \ell e^{-t^2/2}\int_0^{|t|} \frac{u^2}2 \exp\Big\{
%\varkappa u^3\Big(\ell+\frac1{\sqrt n}\Big) +
%\frac{u^2}{2n}\Big(1-2\varkappa u\Big(\ell+\frac1{\sqrt n}\Big)\Big)
%\Big\}du\equiv r_3(t,\ell,n).
%$$
(recall that $\ell_n=\ell+1/\sqrt n$).

Noticing that
$$
\left|K(t)-\frac i{2\pi t}\right|=\frac12(1-t)\sqrt{1+\Big(\cot\pi
t- \frac1{\pi t}\Big)^2},\quad 0\le t\le 1,
$$
we can estimate $\ud$ for any $n\ge2$ and $F$ with a fixed
Lyapunov fraction $\ell$ as
$$
\ud\le 2\int_0^{t_0} |K(t)|\cdot  r_1(Tt,\ell,n)\,dt+
2\int_{t_0}^{1} |K(t)|\cdot f_1(Tt,\ell+1/\sqrt{n},n)\,dt+
$$
$$
+ \frac1{\pi}\int_{t_0}^\infty e^{-T^2t^2/2}\frac{dt}t+
\int_0^{t_0}(1-t)\sqrt{1+\Big(\cot\pi t- \frac1{\pi
t}\Big)^2}e^{-T^2t^2/2}dt\equiv D(\ell,n,t_0,T)
$$
with arbitrary positive $T$ and $t_0\le1$.

The following lemma makes it possible to bound above the set of
the values of $n$ under consideration when estimating the constant
$C_k$ in inequality (3) with $0< k\le 1$.

\smallskip

{\sc Lemma 4}. {\it For any positive $N$, $k\le1$ and
$\eps>k/\sqrt{N}$ for all $t\in\R$ the following estimates hold:}
$$
\sup_{n\ge N}f_j\Big(t,\eps+\frac{1-k}{\sqrt n}\Big)=
f_j\Big(t,\eps+\frac{1-k}{\sqrt N}\Big)\equiv
\widetilde f_{j,N}(t,\eps),\quad j=1,2,
$$
$$
\sup_{n\ge N}r_2\Big(t,\eps-\frac{k}{\sqrt n},n\Big)\le
\eps e^{-t^2/2}\int_0^{|t|} \frac{u^2}2 \exp\Big\{\frac{u^2}{2}-\frac{N-1}{N}\chi\Big(u,\eps+\frac{1-k}{\sqrt N}\Big)\Big\} du
\equiv \widetilde r_{2,N}(t,\eps),
$$
{\it For
$$
|t|\le T(N,\eps)\equiv
\min\Big\{N^{1/4}\eps^{-1/2},(2\varkappa\eps)^{-1}\Big\}
$$
we also have the estimate}
$$
\sup_{n\ge N}r_3\Big(t,\eps-\frac{1}{\sqrt n},n\Big)\le
\frac1{6\varkappa}\Big(e^{\varkappa\eps|t|^3}-1\Big) e^{-t^2/2}
\equiv \widetilde r_{3}(t,\eps).
$$

\smallskip

\proof. The first two statements are trivial consequences of the
monotonicity of the functions $\chi(t,\eps+(1-k)/\sqrt{n})$ and
$f_j(t,\eps+(1-k)/\sqrt{n})$, $j=1,2,$ with respect to $n\ge1$.

To prove the third statement note that the function $r_3$ can be
written in the form
$$
r_3\Big(t,\eps-\frac{1}{\sqrt n},n\Big)=e^{-t^2/2}\int_{0}^{|t|}  \frac{u^2}2 \exp\big\{
\varkappa\eps u^3 + g(n,u)\big\}du,
$$
where
$$
g(x,u)=\ln\Big(\eps-\frac1{\sqrt x}\Big)+
\frac{a(u)}{x},\ x>0,\quad a(u)=\frac{u^2}{2}\Big(1-2\varkappa\eps u\Big),\ u>0.
$$
Since $|t|\le (2\varkappa\eps)^{-1}$ under the conditions of the
lemma, we have $a(u)\ge0$ for all $u\le|t|$. Let us establish that
$g(x,u)$ monotonically increases in $x\ge N$ and $u\le T(N,\eps)$.
Indeed, the derivative
$$
\frac {\partial g(x,u)}{\partial x}=\frac{1}{2x(\eps\sqrt{x}-1)}-
\frac{a(u)}{x^2}
$$
is non-negative if and only if $x-2a(u)\eps\sqrt{x}+2a(u)\ge0$.
Since $a(u)\ge0$, the last condition is satisfied, if
$\sqrt{x}\ge2a(u)\eps\equiv\eps u^2(1-2\varkappa\eps u)$, or,
particulary, if $\sqrt{x}\ge\eps u^2$. So, for all $x\ge N$ and
$u\le T(N,\eps)$ with $T(N,\eps)$ defined in the formulation of
the lemma the function $g(x,u)$ monotonically increases in $x\ge
N$, whence it follows that
$$
\sup_{n\ge N} g(n,u)=\lim_{n\to\infty} g(n,u)=\ln\eps,\quad 0\le u\le T(N,
\eps),
$$
and
$$
\sup_{n\ge N} r_3\Big(t,\eps-\frac{1}{\sqrt n},n\Big)\le \frac12\eps e^{-t^2/2}\int_{0}^{|t|}{u^2} e^{\varkappa\eps u^3}du=
\frac1{6\varkappa}\Big(e^{\varkappa\eps|t|^3}-1\Big) e^{-t^2/2}
\equiv \widetilde r_{3}(t,\eps),
$$
Q. E. D.

\smallskip

Finally, the process of computational optimization can be properly
organized with the help of the following statements.

\smallskip

{\sc Lemma 5} (Bhattacharya and Ranga Rao, 1976).  {\it For any
distribution $F$ with zero mean and unit variance we have}
$$
\rho(F,\Phi)\le \sup_{x>0}\left(\Phi(x)-\frac{x^2}{1+x^2}\right)=
0.54093654\ldots
$$

\smallskip

{\sc Lemma 6}. {\it For any $F\in\F_3$ and $n\ge400$ such that
$\beta^3+1\le0.1\sqrt{n}$ the following estimate takes place$:$}
$$
\ud\le 0.2727\cdot\frac{\beta^3}{\sqrt{n}}+\frac{0.2041}{\sqrt{n}}.
$$
The statement of lemma 6 is a result of the algorithm described in
(Prawitz, 1975b) or (Gaponova and Shevtsova, 2009).

\smallskip

Since the function
$$
g(b)=\frac{0.2727b+0.2041}{b+k},\quad b\ge1,
$$
monotonically increases for $k>0.2041/0.2727=0.74\ldots$ and
monotonically decreases for $0\le k\le0.74$, we have
$$
\sup_{b\ge1}g(b)=
\left\{
  \begin{array}{ll}
    0.2727, & k\ge 0.75,\vspace{2mm} \\
    0.4768/(1+k), & k\le 0.74.
  \end{array}
\right.
$$
Thus, from lemma 6 it follows that for all $n$ and $\beta^3$ such
that ${(\beta^3+k)/\sqrt{n}<0.05(1+k)}$ inequality (3) holds with
$C_k=0.2727$ for $k\ge 0.75$ and with $C_k=0.4768/(1+k)$ for $k\le
0.74$. In particular, for $k=0.425$ we have
$$
\ud\le 0.3346\cdot\frac{\beta^3+0.425}{\sqrt n},\quad \mbox{if}\quad \frac{\beta^3+0.425}{\sqrt n}\le0.07125.
$$

The lemmas presented above give the grounds for restricting the
domain of the values of $\eps=(\beta^3+k)/\sqrt{n}$ by a bounded
interval separated from zero (more details will be given below)
and for looking for the constant $C_k$ in the form
$$
C_k=\max_{\eps}C(\eps),\quad C(\eps)=D(\eps)/\eps,\quad
D(\eps)=\sup\left\{D(\eps,n)\colon n\ge n_*\right\},\eqno(11)
$$
where
$$
D(\eps,n) =\inf_{0<t_0\le1,\, T>0} D(\eps-k/\sqrt{n},n,t_0,T),\eqno(12)
$$
$$
n_*=\max\{1,\,\lceil(1+k)^2/\eps^2\rceil\},
$$
here $\lceil x\rceil$~ is the least integer no less than $x$. The
condition $n\ge n_*$ is a consequence of the inequality
$\beta^3\ge1$. For the estimation of the supremum in $n$ in the
definition of $D(\eps)$, lemma 4 is used for $N$ large enough. The
computation of the maximum in $\eps$ is essentially based on the
property of monotonicity in $\eps$ of all the functions used for
the estimation of $|f_n(t)|$ and $r_n(t)$, and hence, on the
monotonicity of the function $D(\eps)=\eps C(\eps)$. This property
makes it possible to estimate $\max_{\eps}C(\eps)$ using the
values of $C(\eps)$ only in a finite number of points. In
particular, the following statement holds.

\smallskip

{\sc Lemma 7}. {\it For all $\eps_2>\eps_1>0$ the following
inequality is true$:$}
$$
\max_{\eps_1\le\eps\le\eps_2}C(\eps)\le
C(\eps_2)\cdot\frac{\eps_2}{\eps_1}.
$$

\subsubsection{Proof of theorem 1}

Denote
$$
\eps=\ell+\frac{0.425}{\sqrt{n}}=\frac{\beta^3+0.425}{\sqrt{n}}.
$$
Then for $\eps\le0.07$ inequality (6) is a consequence of lemma 6,
and for $\eps\ge 1.62\ge0.541/0.335789$ it follows from lemma 5.
Thus, to compute $C_k$ the maximization with respect to $\eps$ in
(11) is conducted on the interval $0.07\le\eps\le1.62$. To compute
the supremum with respect to $n\ge n_*=\lceil(1.425/\eps)^2\rceil$
we use lemma~4 with $N=600$ for $\eps\le0.1$, $N=300$ for
$0.1<\eps\le0.2$ and $N=100$ for $\eps>0.2$. For the mentioned
values of $\eps$ we have $n_*(0.07)=415$, $n_*(0.1)=204$,
$n_*(0.2)=51$. The maximum with respect to $\eps$ is estimated by
lemma~7 and is attained in the two points: $n=5$,
$\eps\approx0.822$ ($\beta^3\approx1.413$, $t_0\approx 0.385$, $T
= 5.755$) and $n=8$, $\eps\approx0.504$ ($\beta^3=1$, $t_0\approx
0.293$, $T = 8.911$). Both extremal values do not exceed
$0.335789$, whence, theorem 1 is proved.

\subsubsection{Proof of theorem 2}

Denote
$$
\eps=\ell+\frac{1}{\sqrt{n}}=\frac{\beta^3+1}{\sqrt{n}}.
$$
Then for $\eps\le0.1$ inequality (7) is a consequence of lemma 6,
and for $\eps\ge 1.78\ge0.541/0.3051$ it follows from lemma 5.
Thus, to compute $C_k$ the maximization with respect to $\eps$ in
(11) is conducted on the interval $0.1\le\eps\le1.78$. To compute
the supremum with respect to $n\ge n_*=\lceil4/\eps^2\rceil$ we
use the last statement of lemma~4 with $N=200$ and
$T(200,\eps)=\min\{5.04/\eps,3.76/\sqrt{\eps}\}$. It turned out,
that the extremal value is attained at $n\to\infty$ and
$\eps\approx0.985$ ($t_0 = 0.356$, $T =6.147$) and it does not
exceed $0.3051$, Q. E. D.

\section{An improvement of the analog of the Berry--Esseen inequality for Poisson random sums}

\subsection{The history of the problem}

In this section we will use theorem 1 to improve the analog of the
Berry--Esseen inequality for Poisson random sums. Let
$X_1,X_2,...$ be independent identically distributed random
variables with
$$
{\sf E}X_1\equiv\mu,\ \ \ {\sf D}X_1\equiv\sigma^2>0\ \ \
\text{and}\ \ \ {\sf E}|X_1|^3\equiv\beta^3<\infty.\eqno(13)
$$
Let $N_{\lambda}$ be a random variable with the Poisson
distribution with parameter $\lambda>0$. Assume that for any
$\lambda>0$ the random variables $N_{\lambda}$ and $X_1,X_2,...$
are independent. Set
$$
S_{\lambda}=X_1+\ldots+X_{N_{\lambda}}
$$
(for definiteness we assume that $S_{\lambda}=0$ if
$N_{\lambda}=0$). Poisson random sums $S_{\lambda}$ are very
popular mathematical models of many real objects. In particular,
in insurance mathematics $S_{\lambda}$ describes the total claim
size under the classical risk process in the <<dynamical>> case.
Many examples of applied problems from various fields where
Poisson random sums are encountered can be found in, say,
(Gnedenko and Korolev, 1996) or (Bening and Korolev, 2002).

It is easy to see that
$$
{\sf E}S_{\lambda}=\lambda\mu,\ \ \ \ {\sf
D}S_{\lambda}=\lambda(\mu^2+\sigma^2).
$$
The distribution function of the standardized Poisson random sum
$$
\widetilde
S_\la\equiv\frac{S_{\lambda}-\lambda\mu}{\sqrt{\lambda(\mu^2+\sigma^2)}}
$$
will be denoted $F_{\lambda}(x)$.

It is well known that under the conditions on the moments of the
random variable $X_1$ given above, the so-called Berry--Esseen
inequality for Poisson random sums holds: there exists an absolute
positive constant $C<\infty$ such that
$$
\rho(F_{\lambda},\Phi)\equiv\sup_x|F_{\lambda}(x)-\Phi(x)|\le
C\frac{\beta^3}{(\mu^2+\sigma^2)^{3/2}\sqrt{\lambda}}.\eqno(14)
$$
Inequality (14) has rather an interesting history. Apparently, it
was first proved in (Rotar, 1972a) and was published in (Rotar,
1972b) with $C=2.23$ (the dissertation (Rotar, 1972a) was not
published whereas the paper (Rotar, 1972b) does not contain a
proof of this result). Later, with the use of a traditional
technique based on the Esseen smoothing inequality this estimate
was proved in (von Chossy, Raррl, 1983) with $C=2.21$ (the authors
of this paper declared that $C=3$ in the formulation of the
corresponding theorem, which is, of course, true, but actually in
the proof of this theorem they obtained the value $C=2.21$).

In the paper (Michel, 1993) the property of infinite divisibility
of compound Poisson distributions was used to prove that the
constant in (14) is the same as that in the classical
Berry--Esseen inequality. Although Shiganov's estimate $C_0\le
0.7655$ (Shiganov, 1986), had been known by that time (the
original paper by Shiganov had been published in Russian even
earlier, in 1982), Michel used the previous record value due to
Van Beek (Van Beek, 1972) and announced in (Michel, 1993) that
$C\le0.8$ in (14). Being not aware of this paper of Michel, the
authors of the paper (Bening, Korolev and Shorgin, 1997) used an
improved version of the Esseen smoothing inequality and obtained
the estimate $C\le 1.99$. As it has been already noted, the method
of the proof used in (Michel, 1993) is based on the fact that if
for the absolute constant $C_0$ in the classical Berry--Esseen
inequality (1) an estimate $C_0\le M$ is known, then inequality
(14) holds with $C=M$. This circumstance was also noted by the
authors of the paper (Korolev and Shorgin, 1997) in which
independently of the paper (Michel, 1993) the same result was
obtained, but with another currently best estimate $M=0.7655$. As
we noted in Section 1, the best known estimate of the absolute
constant in the classical Berry--Esseen inequality was obtained in
(Tyurin, 2009c), (Tyurin, 2009d): $C_0\le 0.4785$. Therefore,
following the logics of the reasoning used in (Michel, 1993) and
(Korolev and Shorgin, 1997) we can conclude that inequality (14)
holds at least with $C=0.4785$.

In this section we show that actually binding the estimate of the
constant $C$ in (14) to the estimate of the absolute constant
$C_0$ in the classical Berry--Esseen inequality is more loose.
Namely, although the best known upper estimate of $C_0$ is
$M=0.4785$ and moreover, although the unimprovable lower estimate
of $C_0$ is $\approx 0.4097...$, inequality (14) actually holds
with $C=0.3051$. Thus, here we improve the result of (Korolev and
Shevtsova, 2010b) where we proved inequality (14) with $C=0.3450$.

\subsection{Auxiliary results}

The following lemma determines the relation between the
distributions and moments of Poisson random sums and the
distributions and moments of sums of a non-random number of
independent summands. This lemma will be the main tool which we
will use to apply the results known for the classical case, to
Poisson random sums.

Here and in what follows the symbol $\eqd$ will stand for the
coincidence of distributions. Also denote $\nu=\lambda/n$.

\smallskip

{\sc Lemma 7}. {\it The distribution of the Poisson random sum
$S_\la$ coincides with the distribution of the sum of a non-random
number $n$ of independent identically distributed random variables
whatever integer $n\ge1$ is:
$$
X_1+\ldots+X_{N_\la}\ \eqd\ Y_{\nu,1}+\ldots+Y_{\nu,n}
$$
where for each $n$ the random variables
$Y_{\nu,1},\ldots,Y_{\nu,n}$ are independent and identically
distributed. Moreover, if the random variable $X_1$ satisfies
conditions $(13)$, then for the moments of the random variable
$Y_{\nu,1}$ the following relations hold}:
$$
\E Y_{\nu,1} =  \mu\nu,\quad \D Y_{\nu,1} = (\mu^2+\sigma^2)\nu,
$$
$$
\E|Y_{\nu,1}- \mu\nu|^3\le \nu\beta^3(1+40\nu)\ \ \mbox{ for }\ \
n\ge\la.
$$

\smallskip

P r o o f. $\,$ The proof is based on the property of infinite
divisibility of a compound Poisson distribution which implies that
for any integer $n\ge1$ the characteristic function $f_{S_\la}(t)$
of the Poisson random sum $S_\la$ can be represented as
$$
f_{S_\la}(t)= \exp\big\{\la(f(t)-1)\big\}=
\big[\exp\big\{\nu(f(t)-1)\big\}\big]^n \equiv \big[
f_{Y_{\nu,1}}(t) \big]^n,
$$
where $f_{Y_{\nu,1}}$ is the characteristic function of the random
variable $Y_{\nu,1}$. Hence, the distribution of each of the
summands $Y_{\nu,1},\ldots,Y_{\nu,n}$ coincides with the
distribution of the Poisson random sum of the original random
variables:
$$
Y_{\nu,k}\ \eqd\ X_1+\ldots+X_{N_{\nu}},\quad k=1,\ldots,n,
$$
where $N_{\nu}$ is the Poisson-distributed random variable with
parameter $\nu$ independent of the sequence $X_1, X_2,\ldots$
Hence we directly obtain the relations for the first and the
second moments of the random variables $Y_{\nu,1}$ and $X_1$. Let
us prove the relation for the third absolute moments. By the
formula of total probability we have
$$
\E|Y_{\nu,1}- \mu\nu|^{3} \le e^{-\nu} \Big( \nu^{3}| \mu|^{3}
+\nu\E|X_1- \mu\nu|^{3}+\sum_{k=2}^\infty\frac{\nu^k}{k!}
\E|X_1+\ldots+X_k- \mu\nu|^{3} \Big).
$$
Consider the second and the third summands on the right-hand side
separately. For this purpose without loss of generality we will
assume that $n\ge\la$, that is, $\nu\le1$. By virtue of the
Minkowski inequality we have
$$
\left(\E|X_1- \mu\nu|^{3}\right)^{1/3}\le (\bet)^{1/3} +|\mu|\nu=
(\bet)^{1/3} \left(1+ \frac{|\mu|\nu}{(\bet)^{1/3}}\right).
$$
Since $\nu\le1$ and the ratio $|\mu|/(\bet)^{1/3}$ does not exceed
1 by virtue of the Lyapunov inequality, we obtain
$$
\E|X_1- \mu\nu|^{3}\le\bet(1+\nu)^3\le
\bet(1+7\nu).
$$
To estimate the third summand notice that the Lyapunov inequality
yields
$$
\Big|\sum_{i=1}^k x_i\Big|^r\le k^{r-1}\sum_{i=1}^k|x_i|^r, \quad
x_i\in\R,\ i=1,\ldots,k,\ r\ge1,
$$
(see, e. g., (Bhattacharya and Ranga Rao, 1976)). With $r=3$, this
inequality implies
$$
\E|X_1+\ldots+X_k- \mu\nu|^{3}\le \E(|X_1|+\ldots+|X_k|+|
\mu|\nu)^{3}\le
$$
$$
\le (k+1)^2(k\bet+(| \mu|\nu)^{3})\le \bet(k+1)^3
$$
(here we took into account that $|\mu|^{3}\le\bet$ and $\nu\le1$).
Thus,
$$
\E|Y_{\nu,1}- \mu\nu|^{3} \le \nu^{3}| \mu|^{3}+ \nu\E|X_1- \mu\nu|^{3}+
\sum_{k=2}^\infty\frac{\nu^k}{k!} \E|X_1+\ldots+X_k- \mu\nu|^{3}\le
$$
$$
\le\nu\bet\big[1 + (8+K)\nu\big]
$$
where
$$
K=\sum_{k=2}^\infty\frac{(k+1)^3}{k!}=15e-9<32.
$$
The lemma is proved.

\smallskip

{\sc Corollary 1}. {\it Under conditions $(13)$ the distribution
of the standardized Poisson random sum $\widetilde S_\la$
coincides with the distribution of the normalized non-random sum
of $n$ random variables whatever integer $n\ge1$ is$:$
$$
\widetilde S_\la \eqd\ \frac1{\sqrt{n}}\sum_{k=1}^n Z_{\nu,k}
$$
where for each $n$ the random variables
$Z_{\nu,1},\ldots,Z_{\nu,n}$ are independent and identically
distributed. Moreover, these random variables have zero
expectation, unit variance and for all $n\ge\la$ their third
absolute moment satisfies the inequality}
$$
\E|Z_{\nu,1}|^3\le
\frac{\beta^3(1+40\nu)\sqrt{n}}{(\mu^2+\sigma^2)^{3/2}\sqrt{\la}}.\eqno(15)
$$

\smallskip

P r o o f. $\,$ According to lemma 7 for any $n$ we have the
representation
$$
\widetilde S_\la =\frac{S_\la-\la\mu}{\sqrt{\la(\mu^2+\sigma^2)}}\
\eqd\ \frac{Y_{\nu,1}+\ldots+Y_{\nu,n}-n
\mu\nu}{\sqrt{(\mu^2+\sigma^2)n\nu}}\ \equiv\
\frac1{\sqrt{n}}\sum_{k=1}^n Z_{\nu,k},
$$
in which the random variables
$$
Z_{\nu,k}\equiv \frac{Y_{\nu,k}- \mu\nu}{\sqrt\nu} =
\frac{Y_{\nu,k}-\E Y_{\nu,k}}{\sqrt{\D Y_{\nu,k}}}
$$
are independent, identically distributed, have zero expectation
and, unit variance. Moreover, by virtue of the same lemma for all
$n\ge\la$ we have the relation
$$
\E|Z_{\nu,1}|^{3}= \frac{\E|Y_{\nu,1}-\E Y_{\nu,1}|^{3}} {(\D
Y_{\nu,1})^{3/2}} \le\frac{\bet(1+40\nu)}{(
\mu^2+\sigma^2)^{3/2}\nu^{1/2}} =\frac{\bet(1+40\nu)\sqrt{n}}{(
\mu^2+\sigma^2)^{3/2}\sqrt{\lambda}}.
$$
The corollary is proved.

\subsection{Main result}

{\sc Theorem 2}. {\it Under conditions $(13)$ for any $\lambda>0$
we have the inequality}
$$
\rho(F_{\lambda},\Phi)\le
\frac{0.3051\beta^3}{(\mu^2+\sigma^2)^{3/2}\sqrt{\lambda}}.
$$

\smallskip

P r o o f. $\,$ From lemma 7 and corollary 1 it follows that for
any integer $n\ge1$
$$
\rho(F_{\lambda},\Phi)=\sup_x\bigg|{\sf
P}\bigg(\frac1{\sqrt{n}}\sum_{k=1}^n
Z_{\nu,k}<x\bigg)-\Phi(x)\bigg|.
$$
Hence, by theorem 1 for an arbitrary integer $n\ge1$ we have
$$
\rho(F_{\lambda},\Phi)\le 0.3051\frac{{\sf
E}|Z_{\nu,1}|^3}{\sqrt{n}}+\frac{0.3051}{\sqrt{n}}.\eqno(16)
$$
Since $n\ge1$ is arbitrary, we can assume that $n\ge \lambda$,
making it possible to use estimate (15) for the specified $n$ and,
in the continuation of (16), to obtain the inequality
$$
\rho(F_{\lambda},\Phi)\le
0.3051\frac{\beta^3(1+40\la/n)}
{(\mu^2+\sigma^2)^{3/2}\sqrt{\la}}+\frac{0.3051}{\sqrt{n}}.
$$
Since here $n\ge\lambda$ is arbitrary, letting $n\to\infty$ we
finally obtain
$$
\rho(F_{\lambda},\Phi)\le
\lim_{n\to\infty}\bigg[0.3051\frac{\beta^3(1+40\la/n)}
{(\mu^2+\sigma^2)^{3/2}\sqrt{\la}}+\frac{0.3051}{\sqrt{n}}\bigg]=
\frac{0.3051\beta^3}{(\mu^2+\sigma^2)^{3/2}\sqrt{\lambda}},
$$
Q. E. D.

\section{Convergence rate estimates in limit theorems for mixed compound Poisson distributions}

\subsection{Preliminaries}

Let $\Lambda_t$ be a positive random variable whose distribution
depends on some parameter $t>0$. The distribution function of
$\Lambda_t$ will be denoted $G_t(x)$: $G_t(x)={\sf
P}(\Lambda_t<x)$. By a mixed Poisson distribution with a
structural distribution $G_t$ we will mean the distribution of the
random variable $N(t)$ which takes values $k=0,1,...$ with
probabilities
$$
{\sf
P}\big(N(t)=k\big)=\frac{1}{k!}\int\limits_{0}^{\infty}e^{-\lambda}\lambda^kdG_t(\lambda),\
\ \ \ k=0,1,2,...
$$
Some special examples of mixed Poisson distributions are very
well-known. The most well-known and most widely used mixed Poisson
distribution is, of course, the negative binomial distribution
(since it was first used in the form of a mixed Poisson
distribution in (Greenwood and Yule, 1920) to model the
frequencies of accidents). This distribution is generated by the
structural gamma-distribution. Other examples of mixed Poisson
distributions are the Delaporte distribution with the shifted
gamma- structural distribution (Delaporte, 1960), the Sichel
distribution with the generalized inverse Gaussian structural
distribution (Holla, 1967), (Sichel, 1971), Willmot, 1987), The
generalized Waring distribution (Irwin, 1968), (Seal, 1978). The
properties of mixed Poisson distributions are described in detail
in (Grandell, 1997) and (Bening and Korolev, 2002).

Let $X_1,X_2,...$ be independent identically distributed random
variables. Assume that the random variables $N(t),X_1,X_2,...$ are
independent for each $t>0$. Set
$$
S(t)=X_1+\ldots+X_{N(t)}
$$
(for definiteness we assume that if $N(t)=0$, then $S(t)=0$). The
random variable $S(t)$ will be called a mixed Poisson random sum
and its distribution will be called compound mixed Poisson.

In what follows we will assume that the random variables
$X_1,X_2,...$ possess three first moments for which we will use
the same notation as in Section 3 (see (13)). The asymptotic
behavior of the distributions of mixed Poisson random sums $S(t)$
when $N(t)$ infinitely grows in some sense, is principally
different depending on whether $\mu=0$ or not.

The convergence in distribution and in probability will be
respectively denoted by the symbols $\Longrightarrow$ and $\pto$.

First consider the case $\mu=0$. In this case the limit
distributions for standardized mixed Poisson sums are scale
mixtures of normal laws. Without loss of generality, unless
otherwise indicated, we will assume that $\sigma^2=1$.

\smallskip

{\sc Theorem 3} (Korolev, 1996), (Bening and Korolev, 2002). {\it
Assume that $\Lambda_t\pto\infty$ as $t\to\infty$. Then, for a
positive infinitely increasing function $d(t)$ there exists a
distribution function $H(x)$ such that
$$
{\sf P}\bigg(\frac{S(t)}{\sqrt{d(t)}}<x\bigg)\Longrightarrow H(x)\
\ \ \ (t\to\infty)
$$
if and only if there exists a distribution function $G(x)$ such
that for the same function $d(t)$
$$
G_t\big(xd(t)\big)\Longrightarrow G(x) \ \ \ (t\to\infty)\eqno(17)
$$
and}
$$
H(x)=\int\limits_{0}^{\infty}\Phi\big(x/\sqrt{y})dG(y),\ \ \
x\in\r.
$$

\smallskip

Now consider the case $\mu\neq0$. This case is important from the
point of view of insurance applications. Recall that, in general,
${\sf D}X_1=\sigma^2$. Assume that there exist numbers
$\ell\in(0,\infty)$ and $s\in(0,\infty)$ such that
$$
{\sf E}\Lambda_t\equiv \ell t,\ \ \ \ {\sf D}\Lambda_t\equiv
s^2t,\ \ \ t>0.\eqno(18)
$$
Then it is easy to make sure that
$$
{\sf E}S(t)=\mu\ell t,\ \ \ \ {\sf
D}S(t)=[\ell(\mu^2+\sigma^2)+\mu^2s^2]t.
$$
In the book (Bening and Korolev, 2002) a general theorem
presenting necessary and sufficient conditions for the convergence
of compound mixed Poisson distributions was proved. The following
theorem is a particular case of that result.

\smallskip

{\sc Theorem 4} (Bening and Korolev, 2002). {\it Let $\mu\neq0$.
In addition to the conditions on the moments of the structural
random variable $\Lambda_t$ assume that $\Lambda_t\pto\infty$ as
$t\to\infty$. Then, as $t\to\infty$, compound mixed Poisson
distributions converge to the distribution of some random variable
$Z$, that is,
$$
\frac{S(t)-\mu\ell
t}{\sqrt{[\ell(\mu^2+\sigma^2)+\mu^2s^2]t}}\Longrightarrow Z,
$$
if and only if there exists a random variable $V$ such that
$$
\frac{\Lambda_t-\ell t}{s\sqrt t}\Longrightarrow V.
$$
Furthermore,}
$$
{\sf P}(Z<x)={\sf
E}\Phi\bigg(x\sqrt{1+\frac{\mu^2s^2}{(\mu^2+\sigma^2)\ell}}-\frac{\mu
sV}{\sqrt{(\si^2+\mu^2)\ell}}\bigg),\ \ \ x\in\r.
$$

\smallskip

It is easy to see that the limit random variable $Z$ admits the
representation
$$
Z\eqd
\bigg[1+\frac{\mu^2s^2}{(\mu^2+\sigma^2)\ell}\bigg]^{-1/2}\cdot
X+\frac{\mu s}{\sqrt{(\mu^2+\sigma^2)\ell+\mu^2s^2}}\cdot V,
$$
where $X$ is a random variable with the standard normal
distribution independent of $V$.

The basic distinctions of the case $\mu\neq0$ from the case of
compound mixed Poisson distributions with zero expectations
considered above are, first, the necessity of non-trivial
centering and different normalization required for the existence
of non-trivial limit laws and, second, the shape of the limit law
which in this case has the form of a location mixture of normal
laws.

\subsection{Convergence rate estimates in limit theorems for mixed compound Poisson distributions with zero mean}

It is easily seen that the distribution of the mixed Poisson
random sum $S(t)$ can be represented as
$$
{\sf P}(S(t)<x)=\int\limits_{0}^{\infty}{\sf
P}\bigg(\sum_{j=1}^{N_{\lambda}}X_j<x\bigg)dG_t(\lambda),\ \ \
x\in\r.\eqno(19)
$$
Recall that here we assume that
$$
{\sf E}X_1=0,\ \ {\sf E}X_1^2=1,\ \ \beta^3={\sf
E}|X_1|^3<\infty.\eqno(20)
$$
Let $d(t)$, $t>0$, be a positive infinitely increasing function.
In this section we will present some estimates of the rate of
convergence in theorem 3.

For $\lambda>0$ denote
$$
\rho(\lambda)=\sup_x\bigg|{\sf
P}\bigg(\frac{1}{\sqrt{\lambda}}\sum_{j=1}^{N_{\lambda}}X_j<x\bigg)-\Phi(x)\bigg|.
$$
Let $G(x)$ be a distribution function such that $G(0)=0$. If
condition (17) holds, then, according to theorem 3, compound mixed
Poisson distribution of the mixed Poisson sum $S(t)$ normalized by
the square root of the function $d(t)$ converges to the scale
mixture of normal laws in which $G(x)$ is the mixing distribution.
Denote
$$
\Delta_t=\sup_x\bigg|{\sf
P}\bigg(\frac{S(t)}{\sqrt{d(t)}}<x\bigg)-\int\limits_{0}^{\infty}\Phi\bigg(\frac{x}{\sqrt{\lambda}}\bigg)dG(\lambda)\bigg|,\
\ \delta_t=\sup_x\big|G_t\big(d(t)x\big)-G(x)\big|.
$$

\smallskip

{\sc Theorem 5.} {\it Assume that conditions $(20)$ hold. Then for
any $t>0$ we have the estimate}
$$
\Delta_t\le 0.3051\cdot\bet{\sf
E}[\Lambda_t]^{-1/2}+0.5\cdot\delta_t.
$$

\smallskip

P r o o f. $\,$ This statement was first proved in the paper
(Gavrilenko and Korolev, 2006) with a slightly worse constant
(also see (Korolev, Bening and Shorgin, 2007). Here we present a
modified version of the proof. By virtue of representation (19) we
have
$$
\Delta_t=\sup_x\Bigg|\il{0}{\infty}{\sf
P}\Bigg(\sum_{j=1}^{N_{\lambda}}X_j<x\sqrt{d(t)}\Bigg)dG_t(\lambda)-\il{0}{\infty}\Phi\bigg(\frac{x}{\sqrt{\lambda}}\bigg)dG(\lambda)\Bigg|=
$$
$$
=\sup_x\Bigg|\il{0}{\infty}{\sf
P}\Bigg(\frac{1}{\sqrt{\lambda}}\sum_{j=1}^{N_{\lambda}}X_j<x\frac{\sqrt{d(t)}}{\sqrt{\lambda}}\Bigg)dG_t(\lambda)-
\il{0}{\infty}\Phi\bigg(\frac{x}{\sqrt{\lambda}}\bigg)dG(\lambda)\Bigg|=
$$
$$
=\sup_x\Bigg|\il{0}{\infty}{\sf P}\Bigg(\frac{1}{\sqrt{\lambda
d(t)}}\sum_{j=1}^{N_{\lambda
d(t)}}X_j<\frac{x}{\sqrt{\lambda}}\Bigg)dG_t\big(\lambda
d(t)\big)\Bigg)-
\int\limits_{0}^{\infty}\Phi\bigg(\frac{x}{\sqrt{\lambda}}\bigg)dG(\lambda)\Bigg|\le
$$
$$
\le\sup_x\Bigg|\il{0}{\infty}\Bigg[{\sf
P}\Bigg(\frac{1}{\sqrt{\lambda d(t)}}\sum_{j=1}^{N_{\lambda
d(t)}}X_j<\frac{x}{\sqrt{\lambda}}\Bigg)-
\Phi\bigg(\frac{x}{\sqrt{\lambda}}\bigg)\Bigg]dG_t\big(\lambda
d(t)\big)\Bigg|+
$$
$$
+\sup_x\Bigg|\il{0}{\infty}\Phi\bigg(\frac{x}{\sqrt{\lambda}}\bigg)d\big[G_t\big(\lambda
d(t)\big)-G(\lambda)\big]\Bigg|.
$$
Continuing this chain of relations with the use of integration by
parts and theorem 2 we further obtain
$$
\Delta_t\le \il{0}{\infty}\sup_x\Bigg|{\sf
P}\Bigg(\frac{1}{\sqrt{\lambda d(t)}}\sum_{j=1}^{N_{\lambda
d(t)}}X_j<x\Bigg)- \Phi(x)\Bigg|dG_t\big(\lambda d(t)\big)+
$$
$$
+\sup_x\Bigg|\il{0}{\infty}\big[G_t\big(\lambda
d(t)\big)-G(\lambda)\big]d_{\lambda}\Phi\bigg(\frac{x}{\sqrt{\lambda}}\bigg)\Bigg|\le
$$
$$
\le
\il{0}{\infty}\rho(\lambda)dG_t(\lambda)+\sup_{\lambda}\big|G_t\big(\lambda
d(t)\big)-G(\lambda)\big|\cdot\sup_x\Bigg|\il{0}{\infty}d_{\lambda}\Phi\bigg(\frac{x}{\sqrt{\lambda}}\bigg)\Bigg|
\le
$$
$$
\le
0.3051\cdot\bet\il{0}{\infty}\frac{1}{\sqrt{\lambda}}dG_t(\lambda)+0.5\cdot
\sup_{\lambda}\big|G_t\big(\lambda d(t)\big)-G(\lambda)\big|=
$$
$$
=0.3051\cdot\bet{\sf
E}[\Lambda_t]^{-1/2}+0.5\cdot\delta_t,
$$
Q. E. D.

\smallskip

As an example of applications of theorem 5 consider the case where
for each $t>0$ the random variable $\Lambda_t$ has the
gamma-distribution. This case is very important in financial
applications for the asymptotic validation of such popular models
of the evolution of financial indexes as variance-gamma L{\'e}vy
processes (VG-processes) (Madan and Seneta, 1990) or two-sided
gamma-processes (Carr, Madan and Chang, 1998).

As is well known, the density of the gamma-distribution with shape
parameter $r>0$ and scale parameter $\sigma>0$ has the form
$$
g_{r,\si}(x)=\frac{\si^r}{\Gamma(r)}e^{-\si x}x^{r-1},\ \ \ x>0.
$$
Thus, the mixed Poisson distribution with the mixing
gamma-distribution has the characteristic function
$$
\psi(t)=\il{0}{\infty}\exp\{y(e^{it}-1)\}\frac{\si^r}{\Gamma(r)}
e^{-\si y}y^{r-1}dy=
$$
$$
=\frac{\si^r}{\Gamma(r)}\il{0}{\infty}\exp\Big\{-\si y\Big(1+
\frac{1-e^{it}}{\si}\Big)\Big\}y^{r-1}dy=
\Big(1+\frac{1-e^{it}}{\si}\Big)^{-r}.
$$
By the re-parametrization
$$
\si=\frac{p}{1-p}\ \ \ \Big(p=\frac{\si}{1+\si}\Big),\ \ \
p\in(0,1),
$$
we finally obtain
$$
\psi(t)=\Big(\frac{p}{1-(1-p)e^{it}}\Big)^r,\
\ \ t\in\r,
$$
which coincides with the characteristic function of the negative
binomial distribution with parameters $r>0$ and $p\in(0,1)$. So,
in the case under consideration for each $t>0$ the random variable
$N(t)$ has the negative binomial distribution with parameters
$r>0$ and $p\in(0,1)$:
$$
{\sf P}\big(N(t)=n\big)= C_{r+n-1}^{n}p^r(1-p)^{n},\ \ \ \
n=0,1,2,\ldots.\eqno(21)
$$
Here $r>0$ and $p\in(0,1)$ are the parameters and for non-integer
$r$ the quantity $C_{r+n-1}^{n}$ is defined as
$$
C_{r+n-1}^{n}\
=\ \frac{\Gamma(r+n)}{n!\cdot\Gamma(r)}.
$$
In particular, with $r=1$, relation (21) determines the geometric
distribution.

The gamma-distribution function with scale parameter $\si$ ana
shape parameter $r$ will be denoted $G_{r,\si}(x)$. It is easy to
see that
$$
G_{r,\si}(x)\equiv G_{r,1}(\si x).\eqno(22)
$$

The random variable with the distribution function $G_{r,\si}(x)$
will be denoted $U(r,\si)$. It is well known that
$$
{\sf E}U(r,\si)=\frac{r}{\si}.
$$

Fix the parameter $r$ and take $U(r,\si)$ as the random variable
$\Lambda_t$ assuming that $t=\si^{-1}$:
$$
\Lambda_t=U(r,t^{-1}).
$$
As a function $d(t)$ take
$$
d(t)\equiv{\sf E}\Lambda_t={\sf E}U(r,t^{-1}).
$$
Obviously, we have
$$
{\sf E}U(r,t^{-1})=rt.
$$
Then with the account of (22) we have
$$
G_t\big(xd(t)\big)={\sf P}(U(r,t^{-1})<xrt)={\sf
P}(U(r,1)<xr)={\sf P}(U(r,r)<x)=G_{r,r}(x).
$$
Note that the distribution function on the right-hand side of the
latter relation does not depend on $t$. Therefore the choice of
$d(t)$ specified above trivially guarantees the validity of
condition (17) of theorem 3. Moreover, in this case $\delta_t=0$
for all $t>0$.

Now calculate ${\sf E}[\Lambda_t]^{-1/2}$ under the condition
$$
r>\frac{1}{2}.\eqno(23)
$$
We have
$$
{\sf E}[\Lambda_t]^{-1/2}={\sf
E}[U(r,t^{-1})]^{-1/2}=\il{0}{\infty}\frac{e^{-x/t}x^{r-3/2}}{t^r\Gamma(r)}dx=
\frac{\Gamma(r-\frac{1}{2})}{\Gamma(r)\sqrt{t}}.
$$
Thus we obtain the following statement which is actually a
particular case of theorem 5.

\smallskip

{\sc Corollary 2.} {\it Let the random variable $\Lambda_t$ have
the gamma-distribution with shape parameter $r>0$ and scale
parameter $\si=1/t$, $t>0$. Assume that conditions $(20)$ and
$(23)$ hold. Then for each $t>0$ we have}
$$
\sup_x\bigg|{\sf
P}(S(t)<x\sqrt{rt})-\il{0}{\infty}\Phi\bigg(\frac{x}{\sqrt{y}}\bigg)dG_{r,r}(y)\bigg|\le
0.3051\frac{\Gamma(r-\frac{1}{2})}{\Gamma(r)}\cdot\frac{\bet}{\sqrt{t}}.
$$

\smallskip

If $r=1$, then the random variable
$$
N(t)=N_1\big(U(1,t^{-1})\big),\ \ \ t\ge0,
$$
has the geometric distribution with parameter $p=t^{-1}$. As this
is so, the limit (as $t\to\infty$) distribution function of the
standardized geometric sum $S(t)$ is the Laplace distribution with
the density
$$
l(x)=\frac{1}{\sqrt{2}}e^{-\sqrt{2}|x|},\ \ \ \ x\in\r.
$$
The distribution function corresponding to the density $l(x)$ will
be denoted $L(x)$,
$$
L(x)=\begin{cases}{\textstyle\frac12}e^{\sqrt{2}x}, & \text{if
$x<0$,}\vspace{2mm}\cr 1-{\textstyle\frac12}e^{-\sqrt{2}x}, &
\text{if $x\ge0$.}\end{cases}
$$

\smallskip

{\sc Corollary 3.} {\it Let the random variable $\Lambda_t$ have
the exponential distribution with parameter $\si=1/t$, $t>0$.
Assume that conditions $(20)$ hold. Then for each  $t>0$}
$$
\sup_x|{\sf P}(S(t)<x\sqrt{t})-L(x)|\le
0.5408\cdot\frac{\beta_3}{\sqrt{t}}.
$$

\section{Convergence rate estimates in limit theorems for mixed compound Poisson distributions with non-zero mean}

Here we will present some estimates of the rate of convergence in
theorem 4.

\subsection{The case of structural random variables with finite variance}

Under assumptions (18) denote
$$
F_t(x) = \p\left(\frac{S(t)-\mu\ell
t}{\sqrt{[\ell(\mu^2+\sigma^2)+\mu^2s^2]t}}< x\right),
$$
$$
\rho_t=\sup_x \left|F_t(x) - {\sf
E}\Phi\bigg(x\sqrt{1+\frac{\mu^2s^2}{(\mu^2+\sigma^2)\ell}}-\frac{\mu
sV}{\sqrt{(\si^2+\mu^2)\ell}}\bigg) \right|,
$$
$$
G^*(v)=\p \left(V <v\right),\ \ \ \widetilde\delta_t =
\sup_{v}\left|\,G_t\big(vs\sqrt{t}+\ell t\big)-G^*(v)\right|.
$$

\smallskip

{\sc Theorem 6.} {\it Let ${\sf E}X_1=\mu\neq0$, ${\sf
D}X_1=\sigma^2$, ${\sf E}|X_1|^3=\beta^3<\infty$, $\e|V|<\infty$.
Then for any $t>0$ we have
$$\rho_t \leq \widetilde\delta_t
+\frac{1}{\sqrt{t}}\cdot\inf_{\epsilon\in(0,1)}\bigg\{\frac{0.3051\beta^3}
{(\mu^2+\sigma^2)^{3/2}\sqrt{(1-\epsilon)\ell}}+\frac{s}{\ell}\bigg(\frac{{\sf
E}|V|}{\epsilon}+Q(\epsilon)\bigg)\bigg|\bigg\},
$$
where}
$$
Q(\epsilon)=\max\bigg\{\frac{1}{\epsilon},\,\frac{\sqrt{1+\epsilon}}{\big(1+\sqrt{1-\epsilon}\big)\sqrt{2\pi
e(1-\epsilon)}}\bigg\}.
$$

\smallskip

P r o o f. $\,$ A similar statement with slightly worse constants
was first proved in the paper (Artyukhov and Korolev, 2008). Here
we present a modified version of the proof. As above, let
$N_{\lambda}$ be a random variable with the Poisson distribution
with parameter $\lambda>0$ independent of the sequence
$X_1,X_2,\ldots$ Then we can write
$$
\rho_t=\sup_x \left|\p\left(\frac{S(t)-\mu\ell
t}{\sqrt{[\ell(\mu^2+\sigma^2)+\mu^2s^2]t}}<x\right) - \e
\Phi\bigg(x\sqrt{1+\frac{\mu^2s^2}{(\mu^2+\sigma^2)\ell}}-\frac{\mu
sV}{\sqrt{(\mu^2+\sigma^2)\ell}}\bigg) \right|=$$
$$=\sup_x \Bigg|\int\limits_0^\infty\p\left(\frac{S_{N_\lambda} -\mu\ell t}{\sqrt{[\ell(\mu^2+\sigma^2)+\mu^2s^2]t}}
< x\right)dG_t(\lambda)-\hspace{4cm}
$$
$$\hspace{5cm}-\e
\Phi\bigg(x\sqrt{1+\frac{\mu^2s^2}{(\mu^2+\sigma^2)\ell}}-\frac{\mu
sV}{\sqrt{(\mu^2+\sigma^2)\ell}}\bigg) \Bigg|.
$$
Fix an arbitrary $\epsilon\in(0,1)$. Then
$$\rho_t=\sup_x \Bigg|\ \int\limits_{\lambda<(1-\epsilon)\ell t}\p\left(\frac{S_{N_\lambda}
-\mu\ell t}{\sqrt{[\ell(\mu^2+\sigma^2)+\mu^2s^2]t}} <
x\right)dG_t(\lambda)+
$$
$$
+\int\limits_{\lambda>(1+\epsilon)\ell t}\p\left(\frac{S_{N_\lambda}
-\mu\ell t}{\sqrt{[\ell(\mu^2+\sigma^2)+\mu^2s^2]t}} <
x\right)dG_t(\lambda)+
$$
$$
+\int\limits_{(1-\epsilon)\ell t\le\lambda\le(1+\epsilon)\ell t}\p\left(\frac{S_{N_\lambda}
-\mu\ell t}{\sqrt{[\ell(\mu^2+\sigma^2)+\mu^2s^2]t}} <
x\right)dG_t(\lambda)-
$$
$$\hspace{5cm}-\e
\Phi\bigg(x\sqrt{1+\frac{\mu^2s^2}{(\mu^2+\sigma^2)\ell}}-\frac{\mu
sV}{\sqrt{(\mu^2+\sigma^2)\ell}}\bigg) \Bigg|\le
$$
$$
\le \sup_x\Bigg|\ \int\limits_{\lambda<(1-\epsilon)\ell
t}\p\left(\frac{S_{N_\lambda} -\mu\ell
t}{\sqrt{[\ell(\mu^2+\sigma^2)+\mu^2s^2]t}}<
x\right)dG_t(\lambda)+\hspace{4cm}
$$
$$
+\int\limits_{\lambda>(1+\epsilon)\ell t}\p\left(\frac{S_{N_\lambda}
-\mu\ell t}{\sqrt{[\ell(\mu^2+\sigma^2)+\mu^2s^2]t}}<
x\right)dG_t(\lambda)\Bigg|+
$$
$$
+\sup_x\Bigg|\ \int\limits_{(1-\epsilon)\ell
t\le\lambda\le(1+\epsilon)\ell t}\p\left(\frac{S_{N_\lambda} -\mu\ell
t}{\sqrt{[\ell(\mu^2+\sigma^2)+\mu^2s^2]t}}<
x\right)dG_t(\lambda)-
$$
$$
\hspace{5cm}-\e\Phi\bigg(x\sqrt{1+\frac{\mu^2s^2}{(\mu^2+\sigma^2)\ell}}-\frac{\mu
sV}{\sqrt{(\mu^2+\sigma^2)\ell}}\bigg) \Bigg|\le
$$
$$
\le{\sf P}\left(\left|\frac{\Lambda_t}{\ell
t}-1\right|>\epsilon\right)+\hspace{11cm}
$$
$$
+\sup_x\Bigg|\
\int\limits_{(1-\epsilon)\ell t\le\lambda\le(1+\epsilon)\ell
t}\p\left(\frac{S_{N_\lambda} -\mu\ell
t}{\sqrt{[\ell(\mu^2+\sigma^2)+\mu^2s^2]t}}<
x\right)dG_t(\lambda)-
$$
$$
%\begin{equation}
%\label{StInEq0}
\hspace{5cm}-\e\Phi\bigg(x\sqrt{1+\frac{\mu^2s^2}{(\mu^2+\sigma^2)\ell}}-\frac{\mu
sV}{\sqrt{(\mu^2+\sigma^2)\ell}}\bigg) \Bigg|.\eqno(24)
%\end{equation}
$$
Further,
$$ \sup_x\Bigg|\
\int\limits_{(1-\epsilon)\ell t\le\lambda\le(1+\epsilon)\ell
t}\p\left(\frac{S_{N_\lambda} -\mu\ell
t}{\sqrt{[\ell(\mu^2+\sigma^2)+\mu^2s^2]t}}<
x\right)dG_t(\lambda)-
$$
$$\hspace{5cm}-\e\Phi\bigg(x\sqrt{1+\frac{\mu^2s^2}{(\mu^2+\sigma^2)\ell}}-\frac{\mu sV}{\sqrt{(\mu^2+\sigma^2)\ell}}\bigg)
\Bigg|\le
$$
$$
\le\sup_x\ \int\limits_{(1-\epsilon)\ell
t\le\lambda\le(1+\epsilon)\ell t}\Bigg|\
\p\Bigg(\frac{S_{N_\lambda}-\mu\lambda}
{\sqrt{\lambda(\mu^2+\sigma^2)}}<\hspace{6cm} $$
$$\hspace{3cm}<\sqrt{\frac{\ell
t}{\lambda}}\Bigg(x\sqrt{1+\frac{\mu^2s^2}{(\mu^2+\sigma^2)\ell}}-\frac{\mu
s}{\sqrt{(\mu^2+\sigma^2)\ell}} \cdot\frac{\lambda-\ell
t}{s\sqrt{t}}\Bigg)\Bigg)-$$
$$-\Phi\Bigg(\sqrt{\frac{\ell
t}{\lambda}}\Bigg(x\sqrt{1+\frac{\mu^2s^2}{(\mu^2+\sigma^2)\ell}}-\frac{\mu
s}{\sqrt{(\mu^2+\sigma^2)\ell}} \cdot\frac{\lambda-\ell
t}{s\sqrt{t}}\Bigg)\Bigg)\Bigg|\ dG_t(\lambda)+$$
$$+\sup_x\Bigg|\ \int\limits_{(1-\epsilon)\ell t\le\lambda\le(1+\epsilon)\ell t}\!\!\!\!\!\Phi\Bigg(\sqrt{\frac{\ell
t}{\lambda}}\Bigg(x\sqrt{1+\frac{\mu^2s^2}{(\mu^2+\sigma^2)\ell}}-\frac{\mu
s}{\sqrt{(\mu^2+\sigma^2)\ell}} \cdot\frac{\lambda-\ell
t}{s\sqrt{t}}\Bigg)\Bigg)dG_t(\lambda)-
$$
$$-
\e\Phi\bigg(x\sqrt{1+\frac{\mu^2s^2}{(\mu^2+\sigma^2)\ell}}-\frac{\mu
sV}{\sqrt{(\mu^2+\sigma^2)\ell}}\bigg) \Bigg|\equiv I_1+I_2.
$$
Consider $I_1$. Denote
$$
y=\sqrt{\frac{\ell
t}{\lambda}}\Bigg(x\sqrt{1+\frac{\mu^2s^2}{(\mu^2+\sigma^2)\ell}}-\frac{\mu
s}{\sqrt{(\mu^2+\sigma^2)\ell}} \cdot\frac{\lambda-\ell
t}{s\sqrt{t}}\Bigg).
$$
Then $I_1$ can be rewritten in the form
$$
I_1=\sup_y\int\limits_{(1-\epsilon)\ell
t\le\lambda\le(1+\epsilon)\ell t}\Bigg|\,{\sf
P}\left(\frac{S_{N_{\lambda}}-\mu\lambda}{\sqrt{\lambda(\mu^2+\sigma^2)}}<y\right)-\Phi(y)\Bigg|\,dG_t(\lambda)\le$$
$$
\le\int\limits_{(1-\epsilon)\ell t\le\lambda\le(1+\epsilon)\ell
t}\sup_y\Bigg|\,{\sf
P}\left(\frac{S_{N_{\lambda}}-\mu\lambda}{\sqrt{\lambda(\mu^2+\sigma^2)}}<y\right)-\Phi(y)\Bigg|\,dG_t(\lambda).
$$
To estimate the integrand on the right-hand side of the latter
inequality we use theorem 2 and obtain
$$
%\begin{equation}
%\label{I_1}
I_1\le
\frac{0.3051\beta^3}{(\mu^2+\sigma^2)^{3/2}}\int\limits_{\lambda\ge(1-\epsilon)\ell
t}\frac{1}{\sqrt{\lambda}}\,dG_t(\lambda)\le\frac{0.3051\beta^3}{(\mu^2+\sigma^2)^{3/2}\sqrt{(1-\epsilon)\ell
t}}.\eqno(25)
%\end{equation}
$$
Consider $I_2$. We have
$$
I_2\le\sup_x \int\limits_{(1-\epsilon)\ell
t\le\lambda\le(1+\epsilon)\ell t}\Bigg|\,
\Phi\left(\sqrt{\frac{\ell
t}{\lambda}}\Bigg(x\sqrt{1+\frac{\mu^2s^2}{(\mu^2+\sigma^2)\ell}}-\frac{\mu
s}{\sqrt{(\mu^2+\sigma^2)\ell}} \cdot\frac{\lambda-\ell
t}{s\sqrt{t}}\Bigg)\right)-
$$
$$
-\Phi\Bigg(x\sqrt{1+\frac{\mu^2s^2}{(\mu^2+\sigma^2)\ell}}-\frac{\mu
s}{\sqrt{(\mu^2+\sigma^2)\ell}} \cdot\frac{\lambda-\ell
t}{s\sqrt{t}}\Bigg)\Bigg| \ dG_t(\lambda)+
$$
$$
+\sup_x\Bigg|\,\int\limits_{(1-\epsilon)\ell
t\le\lambda\le(1+\epsilon)\ell t}
\Phi\left(x\sqrt{1+\frac{\mu^2s^2}{(\mu^2+\sigma^2)\ell}}-\frac{\mu
s}{\sqrt{(\mu^2+\sigma^2)\ell}} \cdot\frac{\lambda-\ell
t}{s\sqrt{t}}\right)dG_t(\lambda)-
$$
$$
-\e
\Phi\bigg(x\sqrt{1+\frac{\mu^2s^2}{(\mu^2+\sigma^2)\ell}}-\frac{\mu
sV}{\sqrt{(\mu^2+\sigma^2)\ell}}\bigg) \Bigg|\equiv I_{21}+I_{22}.
$$
Denote
$$
z=x\sqrt{1+\frac{\mu^2s^2}{(\mu^2+\sigma^2)\ell}}-\frac{\mu
s}{\sqrt{(\mu^2+\sigma^2)\ell}} \cdot\frac{\lambda-\ell
t}{s\sqrt{t}}.
$$
Then
$$
%\begin{equation} \label{I_21}
I_{21}\le\int\limits_{(1-\epsilon)\ell
t\le\lambda\le(1+\epsilon)\ell
t}\sup_z\bigg|\,\Phi\bigg(z\sqrt{\frac{\ell t}{\lambda}}\bigg)-
\Phi(z)\bigg|\,dG_t(\lambda).\eqno(26)
%\end{equation}
$$
Consider the integrand in (26). By the Lagrange formula we have
$$
%\begin{equation} \label{Lagranj}
\bigg|\,\Phi\bigg(z\sqrt{\frac{\ell t}{\lambda}}\bigg)-
\Phi(z)\bigg|=|z|\cdot\bigg|\sqrt{\frac{\ell
t}{\lambda}}-1\bigg|\,\varphi\bigg(\theta
z+(1-\theta)z\sqrt{\frac{\ell t}{\lambda}}\bigg)\eqno(27)
%\end{equation}
$$
for some $\theta\in[0,1]$ where
$$
\varphi(x)=\Phi'(x)=\frac{1}{\sqrt{2\pi}}e^{-x^2/2}
$$
is the standard normal density. The function
$\varphi(x)=\varphi(|x|)$ monotonically decreases as $|x|$
increases. Therefore the function $\varphi$ on the right-hand side
of (27) attains its maximum value in $\theta\in[0,1]$ at that
value of its argument, whose absolute value is minimum. But the
argument of the function $\varphi$ on the right-hand side of (27)
is itself a linear function of $\theta$. Therefore, the minimum
absolute value of this argument is attained either at $\theta=0$
or at $\theta=1$. But at $\theta=1$ we have
$$
\theta\bigg(1-\sqrt{\frac{\ell t}{\lambda}}\bigg)+\sqrt{\frac{\ell
t}{\lambda}}=1,
$$
while at $\theta=0$ we have
$$
\theta\bigg(1-\sqrt{\frac{\ell t}{\lambda}}\bigg)+\sqrt{\frac{\ell
t}{\lambda}}=\sqrt{\frac{\ell t}{\lambda}}.
$$
In the definition of $I_{21}$ $\lambda$ satisfies the inequality
$\lambda\le(1+\epsilon)\ell t$. Therefore,
$$
\sqrt{\frac{\ell t}{\lambda}}\ge\frac{1}{\sqrt{1+\epsilon}}.
$$
Hence,
$$
\theta\bigg(1-\sqrt{\frac{\ell t}{\lambda}}\bigg)+\sqrt{\frac{\ell
t}{\lambda}}\ge\min\bigg\{1,\,\frac{1}{\sqrt{1+\epsilon}}\bigg\}=
\frac{1}{\sqrt{1+\epsilon}}.
$$
Therefore in $I_{21}$ we have (see (27))
$$
%\begin{equation}
%\label{I_21_1}
\sup_z\bigg|\,\Phi\bigg(z\sqrt{\frac{\ell t}{\lambda}}\bigg)-
\Phi(z)\bigg|\le\bigg|\sqrt{\frac{\ell
t}{\lambda}}-1\bigg|\cdot\sup_z|z|\varphi\bigg(\frac{z}{\sqrt{1+\epsilon}}\bigg).\eqno(28)
%\end{equation}
$$
Furthermore,
$$
\left(z\varphi\bigg(\frac{z}{\sqrt{1+\epsilon}}\bigg)\right)'=
\varphi\bigg(\frac{z}{\sqrt{1+\epsilon}}\bigg)\bigg(1-\frac{z^2}{1+\epsilon}\bigg).
$$
Therefore the supremum in (28) is attained at
$z=\pm\sqrt{1+\epsilon}$ and equals
$$
\sup_z|z|\varphi\bigg(\frac{z}{\sqrt{1+\epsilon}}\bigg)=\sqrt{\frac{1+\epsilon}{2\pi
e}}.
$$
Thus,
$$
I_{21}\le\sqrt{\frac{1+\epsilon}{2\pi
e}}\int\limits_{(1-\epsilon)\ell t\le\lambda\le(1+\epsilon)\ell
t}\bigg|\sqrt{\frac{\ell t}{\lambda}}-1\bigg|\,dG_t(\lambda)=
$$
$$
=\sqrt{\frac{1+\epsilon}{2\pi e}}\int\limits_{(1-\epsilon)\ell
t\le\lambda\le(1+\epsilon)\ell t}\bigg|\frac{\sqrt{\ell
t}-\sqrt{\lambda}}{\sqrt{\lambda}}\bigg|\cdot
\bigg|\frac{\sqrt{\ell t}+\sqrt{\lambda}}{\sqrt{\ell
t}+\sqrt{\lambda}}\bigg|\,dG_t(\lambda)=
$$
$$
=\sqrt{\frac{1+\epsilon}{2\pi e}}\int\limits_{(1-\epsilon)\ell
t\le\lambda\le(1+\epsilon)\ell t}\frac{|\lambda-\ell
t|}{\sqrt{\lambda}\big(\sqrt{\ell
t}+\sqrt{\lambda}\big)}\,dG_t(\lambda)\le
$$
$$
\le\sqrt{\frac{1+\epsilon}{2\pi e}}\cdot
\frac{1}{\sqrt{1-\epsilon}\big(1+\sqrt{1-\epsilon}\big)}\int\limits_{(1-\epsilon)\ell
t\le\lambda\le(1+\epsilon)\ell t} \bigg|\frac{\lambda}{\ell
t}-1\bigg|\,dG_t(\lambda)=
$$
$$
%\begin{equation} \label{I_21_2}
=\sqrt{\frac{1+\epsilon}{2\pi e(1-\epsilon)}}\cdot
\frac{1}{1+\sqrt{1-\epsilon}}\cdot{\sf
E}\bigg|\frac{\Lambda_t}{\ell
t}-1\bigg|\I\bigg(\bigg|\frac{\Lambda_t}{\ell
t}-1\bigg|\le\epsilon\bigg).\eqno(29)
%\end{equation}
$$
Here the symbol $\I(A)$ denotes the indicator of a set $A$.

Consider $I_{22}$. We have
$$
\int\limits_{\lambda:\,|\frac{\lambda}{\ell t}-1|\le\epsilon}
\Phi\left(x\sqrt{1+\frac{\mu^2s^2}{(\mu^2+\sigma^2)\ell}}-\frac{\mu
s}{\sqrt{(\mu^2+\sigma^2)\ell}} \cdot\frac{\lambda-\ell
t}{s\sqrt{t}}\right)\, dG_t(\lambda)=$$
$$
=\int\limits_{\lambda:\,|\frac{\lambda-\ell
t}{s\sqrt{t}}|\le\epsilon\ell\sqrt{t}/s}
\Phi\left(x\sqrt{1+\frac{\mu^2s^2}{(\mu^2+\sigma^2)\ell}}-\frac{\mu
s}{\sqrt{(\mu^2+\sigma^2)\ell}} \cdot\frac{\lambda-\ell
t}{s\sqrt{t}}\right)\, dG_t(\lambda)=$$
$$
=\int\limits_{|v|\le\epsilon\ell\sqrt{t}/s}
\Phi\left(x\sqrt{1+\frac{\mu^2s^2}{(\mu^2+\sigma^2)\ell}}-\frac{\mu
sv}{\sqrt{(\mu^2+\sigma^2)\ell}} \right)\, dG_t(vs\sqrt{t}+\ell
t).
$$
Therefore,
$$%\begin{multline*}
I_{22}=\sup_x\Bigg|\int\limits_{|v|\le\epsilon\ell\sqrt{t}/s}
\Phi\left(x\sqrt{1+\frac{\mu^2s^2}{(\mu^2+\sigma^2)\ell}}-\frac{\mu
sv}{\sqrt{(\mu^2+\sigma^2)\ell}}\right)\, dG_t(vs\sqrt{t}+\ell t)-
$$
$$
-\int\limits_{-\infty}^{\infty}\Phi\left(x\sqrt{1+\frac{\mu^2s^2}{(\mu^2+\sigma^2)\ell}}-
\frac{\mu sv}{\sqrt{(\mu^2+\sigma^2)\ell}}\right)\, dG^*(v)\le
$$%\end{multline*}
$$
\le\sup_x\Bigg|\int\limits_{|v|\le\epsilon\ell\sqrt{t}/s}
\Phi\left(x\sqrt{1+\frac{\mu^2s^2}{(\mu^2+\sigma^2)\ell}}-
\frac{\mu sv}{\sqrt{(\mu^2+\sigma^2)\ell}}\right)
\big[dG_t(vs\sqrt{t}+\ell t)-G^*(v)\big]\Bigg|+$$
$$
+\sup_x\Bigg|\int\limits_{|v|>\epsilon\ell\sqrt{t}/s}
\Phi\left(x\sqrt{1+\frac{\mu^2s^2}{(\mu^2+\sigma^2)\ell}}-
\frac{\mu sv}{\sqrt{(\mu^2+\sigma^2)\ell}}\right)\,
dG^*(v)\Bigg|\equiv I_{221}+I_{222}.
$$
By integration by parts we obtain
$$
I_{221}\le\sup_x\Bigg|\int\limits_{|v|\le\epsilon\ell\sqrt{t}/s}\!\!\!\!\!\!\big[G_t(vs\sqrt{t}+\ell
t)-G^*(v)\big]%\times\hspace{3cm}$$
%$$\hspace{6cm}\times
\,d_v\Phi\left(x\sqrt{1+\frac{\mu^2s^2}{(\mu^2+\sigma^2)\ell}}-
\frac{\mu sv}{\sqrt{(\mu^2+\sigma^2)\ell}}\right)\Bigg|\le
$$
$$
%\begin{equation}
%\label{I_221}
\le\sup_v\big|\,G_t(vs\sqrt{t}+\ell
t)-G^*(v)\big|\equiv\widetilde\delta_t.\eqno(30)
%\end{equation}
$$
Note that in (24) we can apply the Markov inequality and obtain
the estimate
$$
%\begin{equation} \label{StInEq1}
{\sf P}\bigg(\bigg|\frac{\Lambda_t}{\ell
t}-1\bigg|>\epsilon\bigg)\le \frac{1}{\epsilon}\cdot{\sf
E}\bigg|\frac{\Lambda_t}{\ell
t}-1\bigg|\I\bigg(\bigg|\frac{\Lambda_t}{\ell
t}-1\bigg|>\epsilon\bigg).\eqno(31)
%\end{equation}
$$
Further, again applying the Markov inequality we can make sure
that
$$
%\begin{equation} \label{I_222}
I_{222}\le{\sf P}\big(|V|>\epsilon\ell
\sqrt{t}/s\big)\le\frac{s{\sf
E}|V|}{\epsilon\ell\sqrt{t}}.\eqno(32)
%\end{equation}
$$
Now unifying (24), (25), (26), (29), (30), (31) and (32) we
finally obtain that for any $\epsilon\in(0,1)$ there holds the
inequality
$$%\begin{multline*}
\rho_t\le\widetilde\delta_t+\frac{1}{\sqrt{t}}\bigg(\frac{0.3051\beta^3}{(\mu^2+\sigma^2)^{3/2}\sqrt{(1-\epsilon)\ell}}+
\frac{s}{\epsilon\ell}\cdot{\sf E}|V|\bigg) +
$$
$$
+{\sf E}\left|\frac{\Lambda_t}{\ell
t}-1\right|\cdot\max\bigg\{\frac{1}{\epsilon},\,
\sqrt{\frac{1+\epsilon}{2\pi e(1-\epsilon)}}\cdot
\frac{1}{1+\sqrt{1-\epsilon}}\bigg\},
$$
%\end{multline*}
whence we obviously obtain the statement of the theorem since the
Lyapunov inequality obviously implies that for each $t>0$
$$
{\sf E}\left|\frac{\Lambda_t}{\ell
t}-1\right|=\frac{s}{\ell\sqrt{t}}\,{\sf
E}\left|\frac{\Lambda_t-\ell
t}{s\sqrt{t}}\right|\le\frac{s}{\ell\sqrt{t}}\,\sqrt{{\sf
D}\bigg(\frac{\Lambda_t-\ell
t}{s\sqrt{t}}\bigg)}=\frac{s}{\ell\sqrt{t}}\,.
$$
The theorem is proved.

\smallskip

If we additionally assume that the family of random variables
$$
\bigg\{\bigg|\frac{\Lambda_t-\ell t}{s\sqrt{t}}\bigg|\bigg\}_{t>0}
$$
is uniformly integrable, then by the Lyapunov inequality we obtain
the inequality
$$
{\sf E}|V|=\lim_{t\to\infty}{\sf E}\bigg|\frac{\Lambda_t-\ell
t}{s\sqrt{t}}\bigg|\le\lim_{t\to\infty}\sqrt{{\sf
D}\bigg(\frac{\Lambda_t-\ell t}{s\sqrt{t}}\bigg)}=1.\eqno(33)
$$
Hence, from theorem 6 we obtain the following result.

\smallskip

{\sc Corollary 4.} {\it In addition to the conditions of theorem
$6$, let $(33)$ hold. Then for any $t>0$
$$
\rho_t \leq \widetilde\delta_t
+\frac{1}{\sqrt{t}}\cdot\inf_{\epsilon\in(0,1)}\bigg\{\frac{0.3051\beta^3}{(\mu^2+\sigma^2)^{3/2}\sqrt{(1-\epsilon)\ell}}+
\frac{s}{\ell}\Big(\frac{1}{\epsilon}+Q(\epsilon)\Big)\bigg\},
$$
where}
$$
Q(\epsilon)=\max\bigg\{\frac{1}{\epsilon},\,\frac{\sqrt{1+\epsilon}}{\big(1+\sqrt{1-\epsilon}\big)\sqrt{2\pi
e(1-\epsilon)}}\bigg\}.
$$

\subsection{The case of structural random variables with infinite variance}

Assumption (18) which guarantees the existence of the variance of
the structural random variable $\Lambda_t$ is not crucial. An
analog of theorem 6 can be proved for the case where only the
existence of the mathematical expectation of $\Lambda_t$ is
assumed. Namely, the following theorem holds.

\smallskip

{\sc Theorem 7.} {\it Let $\mu\neq0$. Assume that ${\sf
E}\Lambda_t \equiv t$ and $\Lambda_t\pto\infty$ as $t\to\infty$.
Then, as $t\to\infty$, the distributions of normalized mixed
Poisson random sums converge to the distribution of some random
variable $Z$, that is,
$$
\frac{S(t)-\mu}{\sqrt{t}}\Longrightarrow Z,
$$
if and only if there exists a random variable $V$ such that
$$
\frac{\Lambda_t-t}{\sqrt t}\Longrightarrow V.
$$
Moreover,}
$$
{\sf P}(Z<x)={\sf E}\Phi\Bigl(\frac{x-\mu
V}{\sqrt{\si^2+\mu^2}}\Bigr),\ \ \ x\in\r.\eqno(34)
$$

\smallskip

Relation (34) means that in theorem 8
$$
Z\eqd\sqrt{\mu^2+\sigma^2}\cdot X+\mu V
$$
where the random variables $X$ and $V$ are independent and $X$ has
the standard normal distribution.

\smallskip

By analogy with the notation introduced above, denote
$$
\widetilde\rho_t=\sup_x \left|\p\left(\frac{S(t)-\mu t}{\sqrt{t}}<
x\right)-{\sf E}\Phi\bigg(\frac{x-\mu V}{\sqrt{\si^2+\mu^2}}\bigg)
\right|,
$$
$$
\widehat\delta_t =
\sup_{v}\left|\,G_t(v\sqrt{t}+t))-G^*(v)\right|.
$$

\smallskip

{\sc Theorem 8.} {\it Assume that $\beta^3<\infty$, ${\sf
E}\Lambda_t\equiv t$, $t>0$, and $\e|V|<\infty$. Then
$$
\widetilde\rho_t \leq \widehat\delta_t
+\frac{1}{\sqrt{t}}\cdot\inf_{\epsilon\in(0,1)}\bigg\{\frac{0.3051\beta^3}{(\mu^2+\sigma^2)^{3/2}\sqrt{1-\epsilon}}+\frac{{\sf
E}|V|}{\epsilon}+Q(\epsilon){\sf
E}\bigg|\frac{\Lambda_t-t}{\sqrt{t}}\bigg|\bigg\},
$$
where}
$$
Q(\epsilon)=\max\bigg\{\frac{1}{\epsilon},\,\frac{\sqrt{1+\epsilon}}{\big(1+\sqrt{1-\epsilon}\big)\sqrt{2\pi
e(1-\epsilon)}}\bigg\}.
$$

\smallskip

The proof of theorem 8 differs from the proof of theorem 6 only in
notation.

\smallskip

As an example of the situation in which theorems 7 and 8 are
valid, but theorems 4 and 6 are not, consider the case where
$$
\Lambda_t=\max\{0,\sqrt{t}V+t\}+\frac{1}{2t^{\alpha/2}}\Big(\frac{2\alpha+1}{\alpha}\sqrt{t}-1\Big),
$$
with $2<\alpha<3$ and $V$ being the random variable with the
density
$$
p(x)=\frac{\alpha+1}{2(|x|+1)^{\alpha}},\ \ \ \ x\in\r.
$$
It can be easily verified that ${\sf E}\Lambda_t=t$ for any $t>0$,
but the second moment of $\Lambda_t$ is infinite due to that the
second moment of the random variable $V$ does not exist (and
hence, the second moment of the mixed Poisson random sum $S(t)$
with the structural random variable $\Lambda_t$ does not exist).
However, it can be easily seen that
$$
\frac{\Lambda_t-t}{\sqrt{t}}=\max\{-\sqrt{t},\,V\}+
\frac{1}{2t^{(\alpha+1)/2}}\Big(\frac{2\alpha+1}{\alpha}\sqrt{t}-1\Big)
\Longrightarrow V
$$
as $t\to\infty$. This case is an illustrative example of an
interesting and non-trivial fact: unlike the classical summation
theory, for sums with a random number of summands (in particular,
for mixed Poisson random sums) with infinite variances the
existence of non-trivial weak limits is possible under the
normalization of order $t^{1/2}$ which is <<standard>> in the
classical theory only for sums with finite variances.

\smallskip The authors have the pleasure to express their gratitude to
Margarita Gaponova who carried out the supplementary computations
resulting in lemma 6.

\vspace{0.5cm}

\small

\centerline{R E F E R E N C E S}

\begin{enumerate}

\item S. V. Artyukhov and V. Yu. Korolev. Estimates of
the rate of convergence of the distributions of compound doubly
stochastic Poisson processes with non-zero mean to shift mixtures
of normal laws. -- {\it Surveys in Industrial and Applied
Mathematics}, 2008, Vol. 15, No. 6. p. 988--998 (in Russian).

\item V. E. Bening, V. Yu. Korolev and S. Ya. Shorgin.
On approximations to generalized Poisson distribution. -- {\it
Journal of Mathematical Sciences}, 1997, Vol.~83, No. 3, p.
360--367.

\item V. Bening and V. Korolev. {\it Generalized Poisson Models
and their Applications in Insurance and Finance}. VSP, Utrecht,
The Netherlands, 2002.

\item V. Bentkus. {\it On the asymptotical behavior of the constant in the
Berry--Esseen inequality}. Preprint 91 -- 078, Universit{\"a}t
Bielefeld, 1991.

\item V. Bentkus. On the asymptotical behavior of the constant in the
Berry--Esseen inequality. -- {\it Journal of Theoretical
Probability}, 1994, Vol. 2, No, 2, p. 211--224.

\item A. C. Berry. The accuracy of the Gaussian approximation to
the sum of independent variables. -- {\it Transactions of the
American Mathematical Society}, 1941, Vol. 49, p. 122--136.

\item R. N. Bhattacharya and R. Ranga Rao. {\it Normal Approximation
and Asymptotic Expansions}. Wiley, New York, 1976.

\item P. P. Carr, D. B. Madan and
E. C. Chang. The Variance Gamma process and option pricing. --
{\it European Finance Review}, 1998, vol. 2, p. 79--105.

\item P. Delaporte. Un probl\`{e}me de tarification de l'assurance accidents
d'automobile examin\'e par la statistique math\'ematique. -- in:
{\it Trans. 16th Intern. Congress of Actuaries, Brussels}, 1960,
Vol. 2, p. 121--135.

\item J.\,E.\,A.\,Dunnage. MR0397839 (53 \#1695) 60F05. {\it Mathematical Reviews on the Web}, American Mathematical Society, 1977.

\item C.-G. Esseen. On the Liapunoff limit of error in the theory of probability.
-- {\it Ark. Mat. Astron. Fys.}, 1942, Vol. A28, No. 9, p. 1--19.

\item C.-G. Esseen. A moment inequality with an application to the central limit
theorem. -- {\it Skand. Aktuarietidskr.}, 1956, Vol. 39, p.
160--170.

\item M. O. Gaponova and I. G. Shevtsova. Asymptotic estimates of the
absolute constant in the Berry--Esseen inequality for
distributions with infinite third moments. -- {\it Informatics and
Its Applications}, 2009, Vol. 3, No. 4, p.~41-56 (in Russian).

\item S. V. Gavrilenko and V. Yu. Korolev. Convergence
rate estimates for mixed Poisson random sums. -- In: {\it Systems
and Means of Informatics}, Special Issue. Publishing House of the
Institute for Informatics Problems, Russian Academy of Sciences,
Moscow, 2006, p. 248--257 (in Russian).

\item B. V. Gnedenko and V. Yu. Korolev. {\it Random
Summation: Limit Theorems and Applications}. CRC Press, Boca
Raton, 1996.

\item J. Grandell. {\it Mixed Poisson Processes}.
Chapman and Hall, London, 1997.

\item M. S. Holla. On a Poisson-inverse Gaussian distribution. -- {\it
Metrika}, 1967, Vol. 11, p. 115--121.

\item J. O. Irwin. The generalized Waring distribution applied to accident
theory. -- {\it Journal of the Royal Statistical Society, Ser. A},
1968, Vol. 130, p. 205--225.

\item V. Yu. Korolev. A general
theorem on the limit behavior of superpositions of independent
random processes with applications to Cox processes. -- {\it
Journal of Mathematical Sciences}, 1996, Vol. 81, No. 5, p.
2951--2956.

\item V. Yu. Korolev and I. G. Shevtsova. An improvement of the
Berry--Esseen inequality. -- {\it Doklady Academii Nauk}, 2010a, to appear (in Russian). English translation: {\it Doklady
Mathematics}, to appear.

\item V. Yu. Korolev and I. G. Shevtsova. A sharpening of the upper
estimate of the absolute constant in the Berry--Esseen inequality
for mixed Poisson random sums. -- {\it Doklady Academii Nauk}, 2010b, to appear (in Russian). English translation: {\it Doklady
Mathematics}, to appear.

\item V. Yu. Korolev and I. G. Shevtsova. On the upper estimate of
the absolute constant in the Berry--Esseen inequality. -- {\it
Theory of Probability and its Applications}. 2009, Vol. 54, No. 4,
to appear.

\item V. Yu. Korolev and S. Ya. Shorgin. On the absolute constant in the remainder
term estimate in the central limit theorem for Poisson random
sums. -- in: {\it Probabilistic Methods in Discrete Mathematics},
Proceedings of the Fourth International Petrozavodsk Conference.
VSP, Utrecht, 1997, p. 305--308.

\item D. B. Madan and E. Seneta. The variance gamma $($V.G.$)$ model for
share market return. -- {\it Journal of Business}, 1990, Vol. 63,
p. 511--524.

\item R. Michel. On Berry--Esseen results for the compound Poisson
distribution. -- {\it Insurance: Mathematics and Economics}, 1993,
Vol.~13, No. 1, p.~35--37.

\item V. I. Paulauskas. On a smoothing inequality. -- {\it Lithuanian
Mathematical Journal}, 1971, Vol. 11, No. 4, p. 861--866. English
translation by American Mathematical Society in {\it Selected
Translations in Mathematical Statistics and Probability}. 1977,
Vol. 4, p. 7--12.

\item V. V. Petrov. {\it Limit Theorems for Sums of Independent
Random Variables}. Nauka, Moscow, 1987 (in Russian).

\item V. V. Petrov. {\it Limit Theorems of Probability Theory.
Sequences of Independent Random Variables.} Clarendon Press,
Oxford, 1995.

\item H.~Prawitz. Limits for a distribution, if the
characteristic function is given in a finite domain. -- {\it
Scand. AktuarTidskr.}, 1972, p.~138--154.

\item H.~Prawitz. {\it Remainder Term Estimation for Convolution of
Identical Components}. -- Unpublished manuscript of the lecture
given on 16 June, 1972 at the Summer School of the Swedish
Statistical Society, 12--21 June, 1972, L{\"o}ttorp, Sweden.

\item H.~Prawitz. Ungleichungen f\"{u}r den absoluten Betrag einer
charakteristischen Funktion.~-- {\it Skand. AktuarTidskr.}, 1973,
No.~1, pp.~11--16.

\item H.~Prawitz. Weitere Ungleichungen f\"{u}r den absoluten
Betrag einer charakteristischen Funktion.~-- {\it Scand. Actuarial
J.}, 1975, No.~1, pp.~21--29.

\item H.~Prawitz. On the remainder in the central
limit theorem.~I.-- {\it Scand. Actuarial J.}, 1975, No.~3,
p.~145--156.

\item G. V. Rotar. {\it Some Problems of Reserves Planning}. Candidate
Thesis, Central Economical and Mathematical Institute, Moscow,
1972 (in Russian).

\item G. V. Rotar. A problem of storage control. -- {\it Theory of
Probability and its Applications.}, 1972, Vol. 17, No. 3, p.
597--599.

\item H. Seal. {\it Survival Probabilities. The Goal of Risk Theory}. Wiley,
Chichester -- New York -- Brisbane -- Toronto, 1978.

\item I. G. Shevtsova. A refinement of the upper estimate
of the absolute constant in the Berry--Esseen inequality.~-- {\it
Theory of Probability and its Applications}, 2006, Vol.~51, No.~3,
p.~622--626.

\item I. G. Shevtsova. On the absolute constant in the Berry--Esseen inequality.
-- In: {\it The Collection of Papers of Young Scientists of the
Faculty of Computational Mathematics and Cybernetics, Moscow State
University}, Issue 5. Publishing House of the Faculty of
Computational Mathematics and Cybernetics, Moscow State
University, Moscow, 2008, p. 101--110 (in Russian).

\item I. G. Shevtsova. The lower asymptotically exact constant in the central limit theorem.
-- {\it Doklady Academii Nauk}, 2010a, Vol. 430, No.
4, to appear (in Russian). English translation: {\it Doklady
Mathematics}, to appear.

\item I. G. Shevtsova. On a smoothing inequality. -- {\it Doklady Academii Nauk}, 2010b, Vol. 430, No.
5, to appear (in Russian). English translation: {\it Doklady
Mathematics}, to appear.

\item I. S. Shiganov On a refinement of the upper constant in the
remainder term of the central limit theorem. -- In: {\it Stability
Problems for Stochastic Models. Proceedings of the Seminar}.
Publishing House of the Institute for Systems Studies, Moscow,
1982, p. 109--115 (in Russian).

\item I. S. Shiganov. Refinement of the upper bound of the constant in the central limit theorem. --
{\it Journal of Soviet Mathematics}, 1986, Vol. 35, p. 2545--2550.

\item H. S. Sichel. On a family of discrete distributions particular suited
to represent long tailed frequency data. -- in: {\it Proceedings
of the 3rd Symposium on Mathematical Statistics}. Ed. by N. F.
Laubscher. CSIR, Pretoria, 1971, p. 51--97.

\item I. S. Tyurin. On the convergence rate in Lyapunov's
theorem. -- {\it Theory of Probability and its Applications}, to
appear.

\item I. S. Tyurin. On the accuracy of the Gaussian
approximation. -- {\it Doklady Academii Nauk}, 2009, Vol. 429, No.
3, p. 312--316 (in Russian). English translation: {\it Doklady
Mathematics}, 2009, Vol. 80, No. 3, to appear.

\item I. S. Tyurin. Refinement of the upper bounds of the
constants in Lyapunov's theorem. -- {\it Russian Mathematical
Surveys}, to appear.

\item I. Tyurin. New estimates of the convergence rate in the
Lyapunov theorem. -- arXiv:0912.0726v1, 3 December, 2009.

\item J. D. Vaaler. Some extremal functions in Fourier analysis. -- {\it Bulletin
of the American Mathematical Society (New Series)}, 1985, Vol. 12,
No. 2, p. 183--216.

\item P. van Beek. An application of Fourier methods to the problem of sharpening
the Berry--Esseen inequality. -- {\it Z. Wahrsch. verw. Geb.},
1972, Bd. 23, S. 187--196.

\item R. von Chossy and G. Rappl. Some approximation methods for the
distribution of random sums. -- {\it Insurance: Mathematics and
Economics}. 1983, Vol. 2, p. 251--270.

\item G. E. Willmot. The Poisson--inverse Gaussian distribution as an alternative
to the negative binomial. -- {\it Scandinavian Actuarial Journal},
1987, p. 113-127.

\item V. M. Zolotarev. On closeness of the distributions of two sums of independent random
variables. -- {\it Theory of Probability and its Applications},
1965, Vol. 10, No. 3, p. 519--526.

\item V. M. Zolotarev. An absolute estimate of the remainder term in the central limit theorem,
-- {\it Theory of Probability and its Applications}, 1966, Vol.
11, No. 1, p. 108--119.

\item V. M. Zolotarev. A sharpening of the inequality of Berry--Esseen. --
{\it Z. Wahrsch. verw. Geb.}, 1967, Bd. 8, S. 332--342.

\item V. M. Zolotarev. Some inequalities in probability theory and
their application in sharpening the Lyapunov theorem. --
{\it Soviet Math. Dokl.}, 1967, Vol. 8, p. 1427--1430.

\item V. M. Zolotarev. {\it Modern Theory of Summation of Random
Variables}. VSP, Utrecht, The Netherlands, 1997.

\end{enumerate}

\end{document}